\documentclass[,preprint]{imsart}

\usepackage[a4paper,width=140mm,top=25mm,bottom=25mm]{geometry}

\usepackage{color}
\usepackage{latexsym}
\usepackage{mathtools}
\usepackage{amsmath}
\usepackage{amssymb}
\usepackage[utf8]{inputenc}
\usepackage{textcomp}
\usepackage{graphicx}
\usepackage{epsfig}
\usepackage{hyperref}
\usepackage{bm} 

\usepackage{enumitem}
\usepackage{here}

\usepackage{multirow}
\usepackage{color}
\usepackage[dvipsnames]{xcolor}

\RequirePackage{amsthm,amsmath,amsfonts,amssymb}
\RequirePackage[authoryear]{natbib}
\RequirePackage[colorlinks,citecolor=blue,urlcolor=blue]{hyperref}
\RequirePackage{graphicx}

\startlocaldefs

\providecommand{\skakko}[1]{\left(#1\right)}
\providecommand{\mkakko}[1]{\left\{#1\right\}}

\theoremstyle{plain}
\newtheorem{theorem}{Theorem}[section]

\theoremstyle{plain}
\newtheorem{lem}{Lemma}[section]

\theoremstyle{remark}

\newtheorem{ex}{Example}[section]
\newtheorem{assumption}{Assumption}[section]
\newtheorem{rem}{Remark}[section]

\newcommand{\abs}[1]{\left\lvert#1\right\rvert}
\usepackage{bigints}

\DeclareMathOperator{\argmax}{\textrm{argmax}}

\endlocaldefs

\begin{document}

\begin{frontmatter}
	\title{Residual spectrum: Brain functional connectivity detection beyond coherence}
	\runtitle{Residual spectrum beyond coherence}
	
	\begin{aug} 
\author[A,e1]{\fnms{Yuichi} \snm{Goto}\ead[label=e1,mark]{yuichi.goto@math.kyushu-u.ac.jp}},
\author[B,e2]{\fnms{Xuze} \snm{Zhang}\ead[label=e2,mark]{xzhang51@umd.edu}}
\author[B,e3]{\fnms{Benjamin} \snm{Kedem}\ead[label=e3,mark]{bnk@umd.edu}}
\and 
\author[D,e4]{\fnms{Shuo} \snm{Chen}\ead[label=e4,mark]{shuochen@umd.edu}}

\address[A]{Faculty of Mathematics, Kyushu University  \printead{e1}}
 
\address[B]{Department of Mathematics and Institute for Systems Research, University of Maryland \printead{e2}, \printead{e3}}

\address[D]{Maryland Psychiatric Research Center, School of Medicine, University of Maryland \printead{e4}}

\end{aug}

\begin{abstract}
Coherence is a widely used measure to assess linear relationships between time series. However, it fails to capture nonlinear dependencies. To overcome this limitation, this paper introduces the notion of residual spectral density as a higher-order extension of the squared coherence. The method is based on an orthogonal decomposition of time series regression models. We propose a test for the existence of the residual spectrum and derive its fundamental properties. A numerical study illustrates finite sample performance of the proposed method. An application of the method shows that the residual spectrum can effectively detect brain connectivity. 
Our study reveals a noteworthy contrast in connectivity patterns between schizophrenia patients and healthy individuals. Specifically, we observed that non-linear connectivity in schizophrenia patients surpasses that of healthy individuals, which stands in stark contrast to the established understanding that linear connectivity tends to be higher in healthy individuals. This finding sheds new light on the intricate dynamics of brain connectivity in schizophrenia.
\end{abstract}



\end{frontmatter}

\section{Introduction}

Functional connectivity within the human brain is important for the understanding of brain disorders such as depression, Alzheimer's disease, autism, dyslexia, dementia, epilepsy, and other types of mental disorder. The complete map of the human brain functional connectivity, or neural pathways, is referred to as the {\em connectome}. 
In describing the topology of the human brain, functional connectivity is commonly calculated by characterizing the {\em co-variation} of blood oxygen level-dependent (BOLD) signals, or time series, between two distinct neural populations or regions. To investigate the functional architecture of the brain, time series from resting-state functional magnetic resonance imaging (fMRI) which measures spontaneous low-frequency fluctuations in BOLD signals are often used.
 

Various time and spectral domain methods have been used in the analysis of BOLD signals. Among them is the wide spread measure of coherence  (see, e.g., \cite{mlbv01} and \cite{wwywyd15}). Coherence is a time series index used to determine if two or more brain regions are spectrally connected, that is, if two or more brain regions have ``similar neuronal oscillatory activity'' \citep{bowyer16}. See \citet[Section 11, p.436]{brockwellanddavis1991} for the definition.
However, coherence measures the strength of linear relationships only, but fails to measure nonlinear relationships, especially, the spectral contributions of quadratic interaction terms. The reasons why nonlinearities should be incorporated into the analysis of brain connections are explained in \citet[p.1572]{zhaoetal19} and \citet[Section 5.3]{caoetal22}.

The drawback is rectified by the so called {\em maximum residual coherence}, which provides the contribution of specific interaction terms to co-variation of brain regions, over and beyond the coherence by itself. 
Its origin can be traced back to the orthogonal decomposition of linear and quadratic functionals in the spectral domain proposed by \cite{kimelfeld74}. Based on the decomposition, \cite{kedem75} developed a selection criterion, {\em lagged coherence}, to determine the inclusion of {\em lagged processes} or interaction terms in a linear system. This criterion is a function of frequency so that ``we may run into a situation where one lag maximizes the coherence over a certain
	frequency band while another lag maximizes it over a different band'' \citep{kedem75}. The maximum residual coherence was developed to resolve such issue and shown to be effective in both cases of continuous  and binary  time series (\citealt{kkk14} and \citealt{kedem16}). The criterion is defined using supremum with respect to frequency, while \cite{zk21} proposed an alternative criterion using integration with respect to frequency to choose interaction terms in a more general setting.

Although the method based on maximum residual coherence is promising, {no theoretical framework for its use in testing problems exists so far.} Therefore, in this paper, {we propose a new spectrum, referred to as the {\it residual spectrum} to capture non-linear association}, and fill the gap by proposing, as a special case, a test for the existence of the interaction terms {in covariation.} We demonstrate that 
{the disparity in the test results between between healthy individual and schizophrenia patients of our residual spectrum exceeds than that of the coherence. Additionally, we found that the non-linear connectivity for schizophrenia patients is higher than that for healthy individuals, whereas it is known that the reverse is true for linear connectivity.}

 In addition to coherence, \cite{msnp20} listed alternative measures that can capture co-variation between brain signals in different ways including wavelet coherence, mutual information, dynamic time wraping, and more. Regarding other extensions of coherence, \cite{ov08} proposed the local band coherence to deal with non-stationary signals. \cite{fo11} studied brain functional connectivity by using the partial coherence estimated by the generalized shrinkage estimator, which is based on both nonparametric and parametric spectral density estimation.
 \cite{lenart11} dealt with the magnitude of coherence for almost periodically correlated time series and \cite{bk19} proposed the quantile coherence based on quantile spectra. \cite{eso19} proposed a cluster coherence to measure similarity between vector time series and applied to electroencephalogram (EEG) signals. {See also \cite{Matsuda06} for a test to construct graphical models in the frequency domain and \cite{LTM21} for the concept of local Granger causality and its associated testing procedures.}

However, our goal is to highlight a way to extend and enhance the measure of coherence in the analysis of BOLD signals by {incorporating interaction terms}.  \cite{caoetal22} gave a {comprehensive} review of recent developments on statistical methods to analyze brain functional connectivity.

For hypothesis testing problems described by spectra, $L_2$- and supremum-type statistics have been studied by many authors. 
As for $L_2$-type statistics,  \cite{tpk96} discussed the test based on the integrated function of spectral density matrices (see also \citealt[Chapter 6]{tk00}), which can be applied to, e.g., the test for the magnitude of linear dependence and discriminant analysis (\cite{kst98}). A closely related problem was examined in \cite{ym09} under Gaussian assumption. 
To include other hypotheses, \cite{eichler08} considered the testing problem described by the integrated squared (Euclidean) norm of the function of spectral density matrices. In the same spirit of \cite{eichler08}, tests for the equality of spectra are proposed using the bootstrap method \citep{dp09} and the randomization method \citep{jp15}.  On the other hand, supremum-type statistics are considered by \cite{wv67}, \citet[Thorem 12, Chapter V.5]{hannan1970}, \cite{rudzkis93}, and \cite{wz18}. However, it is difficult to obtain asymptotics of supremum statistics for the function of a spectral density matrix in general, but it can for $L_2$-statistics. Therefore,  in this article, we propose {an} $L_2$-type statistic.

The paper is organized as follows.
{Section \ref{sec2} introduce the residual spectrum, which is the main objective in this paper.}
Section \ref{sec3} {delves into} time series multiple regression models and assumptions, and show there exists a unique orthogonal representation of the model {to derive the residual spectrum and provide an interpretation of it}. 
In Section \ref{sec4}, we propose a test for the existence of the residual spectrum. The asymptotic null distribution of our test statistic and consistency of the test are derived. The consistency of estimators of unknown parameters in the model is also addressed. 
Section \ref{sec5} provides the finite sample performance of our method. Section \ref{sec6} shows the utility of our method by applying the proposed test to brain data. All proofs of Lemma and Theorems in the main article are presented in Appendix \ref{appen}.

{\section{Residual spectrum}\label{sec2}
In this section, we introduce a new spectrum, termed the residual spectrum of order $j$, for the vector stationary process $(X_0(t),X_1(t),\ldots,X_K(t))$,  designed to capture both linear and non-linear terms of $(X_1(t),\ldots,X_K(t)$ in explaining $X_0(t)$.
For the sake of brevity in notation, denote the spectral density matrix for $(X_0(t),X_1(t),\ldots,X_K(t))$ and $(X_1(t),\ldots,X_K(t))$ by
$\bm{f}(\lambda):= \skakko{f_{ij}(\lambda)}_{i,j=0,\ldots,K}$ and $\bm{f}_K(\lambda):=\skakko{f_{ij}(\lambda)}_{i,j=1,\ldots,K}$, respectively.}

{Under the regularity conditions, as derived in the subsequent section, the {\it residual spectrum $f_{G_jG_j}$ of order $j$} is defined as follows: for $j=1$,
\begin{align*}
f_{G_1G_1}(\lambda):=
\frac{\left|f_{10}(\lambda)\right|^2}{f_{11}(\lambda)}\quad \text{and }
f_{G_jG_j}(\lambda)=
\frac{\left|A_{jj}  \skakko{e^{i\lambda}}\right|^2{\rm det}\skakko{\bm{f}_{j}(\lambda)}}{{\rm det}\skakko{\bm{f}_{j-1}(\lambda)}} \quad\text{for $j=2,\ldots,K$,}
\end{align*}
where, for $j\in\{2,\ldots,K\}$, 
\begin{align*}
A_{jj}  \skakko{e^{i\lambda}}
:=
\frac{-\sum_{i=1}^{j-1}{\rm det}\skakko{\overline{\bm{f}_{i,j}^\ddag(\lambda)}}f_{i0}(\lambda)+{\rm det}\skakko{\bm{f}_{j-1}(\lambda)}f_{j0}(\lambda)}{{\rm det}\skakko{\bm{f}_{j}(\lambda)}}
\end{align*}
and 
$$
{\bm{f}}_{i,j}^\ddag(\lambda):=
\begin{pmatrix}
{f}_{11}(\lambda)&\cdots&{f}_{1(i-1)}&{f}_{1j}&{f}_{1(i+1)}&\cdots&{f}_{1(j-1)}\\
\vdots&&\vdots&\vdots&\vdots&&\vdots\\
{f}_{(j-1)1}(\lambda)&\cdots&{f}_{(j-1)(i-1)}&{f}_{(j-1)j}&{f}_{(j-1)(i+1)}&\cdots&{f}_{(j-1)(j-1)}\\
\end{pmatrix}.
$$}

{The interpretation of our spectrum, as elucidated in the subsequent section, is as follows:
\begin{enumerate}
\item
$f_{G_1G_1}(\lambda)$ quantifies the strength of the linear association between $X_0(t)$ and $X_1(t)$. 
\item
$f_{G_2G_2}(\lambda)$ measures the strength of the linear relationship between $X_0(t)$ and $X_2(t)$ after eliminating the linear effect of $X_1(t)$.
\item
$f_{G_3G_3}(\lambda)$ indicates the strength of the linear association between $X_0(t)$ and $X_3(t)$ after removing the linear effects of $X_1(t)$ and $X_2(t)$.
\item Similarly, $f_{G_jG_j}(\lambda)$ represents the strength of the linear relationship between $X_0(t)$ and $X_j(t)$ subsequent to the removal of the linear effects of $X_1(t),\ldots,X_{j-1}(t)$.
\end{enumerate}}

{While one might initially presume that $f_{G_jG_j}(\lambda)$ is incapable of capturing non-linear relationships, it does indeed possess this capability.
We shall elaborate on the integration of non-linear relationships through our residual spectrum. 
Let $X_0(t)$ and $X_1(t)$ be stationary processes and consider taking $X_2(t)$ as a non-linear process of $X_1(t)$, for instance, $X_2(t)=X_1(t)X_1(t+u)$ for some predetermined time lag $u$.
In this scenario, $f_{G_2G_2}(\lambda)$ accounts for the contribution of the non-linear terms $X_2(t)$ which cannot be  measured by the linear term $X_1(t)$. The quantity $f_{G_1G_1}(\lambda)+f_{G_2G_2}(\lambda)$ can be interpreted as the spectral density, encompassing  not solely the linear association of $X_1(t)$ but also the non-linear association of $X_2(t)$. 
This extension can be generalized to higher-order scenarios in a similar fashion. In the following section, we shall provide a derivation of our residual spectrum and elucidate the rationale behind its interpretation as aforementioned.}

\section{Orthogonal decomposition}\label{sec3}
We consider the model
\begin{align}\label{model_K}
X_0(t):=\zeta+\sum_{i=1}^K\sum_{k_i=-\infty}^\infty b_i(k_i)X_i({t-k_i})+
\epsilon(t),
\end{align}
where $X_0(t)$ is a response variable, $\zeta$ is intercept, $\bm{X}(t):=(X_1(t),\ldots,X_K(t))^\top$ is a covariate process with autocovariance matrix $\Gamma_{\bm
X}(u):={\rm E}\bm{X}(t)\bm{X}({t-u})$ such that, for any $i=1,\ldots,K$, $X_i(t):=X_i^\prime-{\rm E}X_i^\prime$ where $(X_1^\prime(t),\ldots,X_K^\prime(t))^\top$ a $s$-th order stationary process for any $s\in\mathbb N$, $\epsilon(t)$ is an i.i.d.\ centered disturbance process independent of $\{\bm{X}{(t)}; t\in \mathbb Z\}$ with moments of all orders, and $\sum_{k=-\infty}^\infty k^2|b_i(k)|<\infty$ for all $i=1,\ldots,K$. 
{The continuous version of this model was considered by \citet[Section 11.4.1, p.485]{wj68}.}

{For any random vector $\bm W(t):=(W_1(t),\cdots,W_\ell(t))^\top$ and $t_1,\ldots,t_\ell\in\mathbb Z$, {the cumulant of order $\ell$} of $(W_1(t_1),\cdots,W_\ell(t_\ell))$ is defined as
\begin{align*}
{\rm cum}(W_1(t_1),\cdots,W_\ell(t_\ell)):=&\\
\sum_{(\nu_1,\ldots,\nu_p)}(-1)^{p-1}(p-1)!&\skakko{{\rm E}\prod_{j\in\nu_1}W_{\nu_1}(t_{\nu_1})}\ldots\skakko{{\rm E}\prod_{j\in\nu_p}W_{\nu_p}(t_{\nu_p})},
\end{align*}
where the summation $\sum_{(\nu_1,\ldots,\nu_p)}$ extends over all partitions $(\nu_1,\ldots,\nu_p)$ of $\{1,2,\cdots,\ell\}$ (see \citealp[p.19]{brillinger1981}).} In this paper, we make the following assumption.
\begin{assumption}\label{as_moment} 
For any integer $\ell\geq 2$ and {$(i_1,\ldots,i_\ell)\in\{1,\ldots,K\}^\ell$}, it holds that
\begin{equation}\label{as_eq_moment}
\sum_{s_2,\ldots,s_\ell=-\infty}^\infty\skakko{1+\sum_{j=2}^\ell\abs{s_j}^2}\abs{\kappa_{i_1,\ldots,i_\ell}(s_2,\ldots,s_\ell)}<\infty,
\end{equation}
where $\kappa_{{i_1},\ldots,{i_\ell}}(s_2,\ldots,s_\ell)
={\rm cum}\{X_{i_1}(0),X_{i_2}(s_2),\ldots,X_{i_\ell}(s_\ell)\}$.
\end{assumption}
In conjunction with $\sum_{k=-\infty}^\infty k^2|b_i(k)|<\infty$ for all $i=1\ldots,K$ and the existence of moments of all orders for $\epsilon(t)$, Assumption \ref{as_moment} implies the cumulant summability \eqref{as_eq_moment} for any integer $\ell\geq 2$ and $(i_1,\ldots,i_\ell)\in\{0,1,\ldots,K\}^\ell$. 
The condition about the summability of cumulant is standard in this context, e.g., \citet[Assumption 2.6.2]{brillinger1981}, \citet[Assumption 2]{tpk96}, \citet[Assumption 3.1]{eichler08}, \citet[Section 4]{shao10}, \citet[Assumption 2.1]{jp15}, and \citet[Section 4]{av20}. The sufficient condition for the summability of the cumulant is often discussed. The summability of the fourth-order cumulant under $\alpha$-mixing condition was shown by \citet[Lemma 1]{andrews91}. For univariate alpha-mixing processes, \citet[Remark 3.1]{neumann96}, \citet[Remark 3.1]{ls17}, and \citet[Lemma 2]{bdh20} showed the summability condition holds for higher-order cases. The essential tool of the proof is Theorem 3 of \cite{sj1988}, which provides the upper bound of higher-order cumulant in terms of alpha-mixing coefficients and moments for univariate processes. Thus, the theorem is not applicable for multivariate processes. Besides, \citet[Proposition 4.1]{pt13} showed the MA($\infty$) process meets the summability condition. On the other hand, \citet[Lemmas 4.1 and 4.2]{kley14} gave a sufficient condition for the summability of cumulant for indicator functions of the processes in a beautiful way (see also \citet[Proposition 3.1]{dette16}). Employing his idea, we obtain the following lemma.
\begin{lem}\label{summability_mixing}
For the geometrically $\alpha$-mixing strictly stationary vector process $\bm{X}(t):=(X_1(t),\ldots,X_K(t))^\top$ with moments of all orders in the sense that $$\sup_{s_2,s_3,\ldots,s_\ell\in\mathbb Z}\abs{{\rm E}X_{i_1}(0)X_{i_2}(s_2)\cdots X_{i_\ell}(s_\ell)}<\infty$$ for any $\ell\in\mathbb N$ and any $(i_1,\ldots,i_\ell)\in\{1,\ldots,K\}^\ell$ with $\alpha$-mixing coefficient $\alpha(\cdot)$ such that $\alpha(n)\leq C_\alpha\rho^n$, 
where
\begin{align*}
\alpha(n):=\sup_{k\in\mathbb Z, A\in \mathcal F_{-\infty}^k,\ B \in \mathcal F_{k+n}^{\infty}}|{\rm P}(AB)-{\rm P}(A){\rm P}(B)|,
\end{align*}
for $a\leq b$, $\mathcal F_{a}^b$ is the $\sigma$-field generated by $\{\bm X(t): a\leq t\leq b\}$, and some constants $C_\alpha\in(1,\infty)$ and $\rho\in(0,1)$, it holds, for any $d\in\mathbb N$, any $\ell\in\mathbb N$, and any $(i_1,\ldots,i_\ell)\in\{1,\ldots,K\}^\ell$,
\begin{align}\label{lem_summability_mixing_eq}
\sum_{s_2,\ldots,s_\ell=-\infty}^\infty\skakko{1+\sum_{j=2}^\ell\abs{s_j}^d}\abs{{\rm cum}\{X_{i_1}(0),X_{i_2}(s_2),\ldots, X_{i_\ell}(s_\ell)\}}<\infty.
\end{align}
\end{lem}
Therefore, the geometrically $\alpha$-mixing strictly stationary vector process with the appropriate moment condition holds {under} Assumption \ref{as_moment}. We present two examples.

\begin{ex}[vector ARMA process]\label{ARMA}{\rm
First example satisfying Assumption \ref{as_moment} is
the causal vector ARMA ($p$,$q$) process $\bm{X}(t):=(X_1(t),\ldots,X_K(t))^\top$, defined as follows: for $p,q\in\mathbb N$, 
\begin{align}\label{ex_ARMA}
\bm\Phi(B)\bm{X}(t)=\bm\Theta(B)\bm{\xi}(t),
\end{align}
where $\bm\Phi(z):=\bm{I}_K-\sum_{k=1}^p\bm\Phi_k z^k$, $\bm\Theta(z):=\bm{I}_K+\sum_{k=1}^q\bm\Theta_kz^k$, $\bm\Phi_k$ and $\bm\Theta_k$ are $K$-by-$K$ matrices satisfying ${\rm det}(\bm\Phi(z))\neq0$ for all $|z|\leq 1$, $B$ is the backward shift operator, $\bm{\xi}(t):=(\xi_1(t),\ldots,\xi_K(t))^\top$ is i.i.d.\ random variables with moments of all orders. Then, $\bm{X}(t)$ satisfies Assumption \ref{as_moment} since  ARMA process has the geometrically $\alpha$-mixing property \citep{mokkadem88} and innovation process satisfies the moment condition stated in Lemma \ref{summability_mixing}. 
}\end{ex}

\begin{ex}[Volterra series]{\rm
Second example is the non-linear time series defined, 
for the geometrically $\alpha$-mixing strictly stationary $K^\prime$-dimensional process $\bm{\upsilon}(t):=(\upsilon_1(t),\ldots,\upsilon_{K^\prime}(t))^\top$ with the moment condition stated in Lemma \ref{summability_mixing} and for $i_{j,d,s}\in\{1,\ldots,K^\prime\}$ $(j\in\{1,\ldots,K\}, d\in\{1,\ldots,D\}, s\in\{1,\ldots,d\})$, as $$\bm{X}(t):=(X_1(t),\ldots,X_K(t))^\top$$ and
\begin{align*}
X_j(t)
:=&
\sum_{k_1=-J}^J c_{j,1}(k_1)\upsilon_{i_{j,1,1}}(t-k_1)
+\sum_{k_1,k_2=-J}^J c_{j,2}(k_1,k_2)\upsilon_{i_{j,2,1}}(t-k_1) \upsilon_{i_{j,2,2}}(t-k_2)\\
&
+\sum_{k_1,k_2,k_3=-J}^J c_{j,3}(k_1,k_2,k_3)\upsilon_{i_{j,3,1}}(t-k_1) \upsilon_{i_{j,3,2}}(t-k_2)\upsilon_{i_{j,3,3}}(t-k_3)+\cdots\\
&
+\sum_{k_1,\ldots,k_j=-J}^J c_{j,j}(k_1,\ldots,k_D)\upsilon_{i_{j,j,1}}(t-k_1)\cdots \upsilon_{i_{j,j,j}}(t-k_j)\\
=&
\sum_{d=1}^j \sum_{k_1,\ldots,k_d=-J}^J c_{j,d}(k_1,\ldots,k_d)\prod_{s=1}^d\upsilon_{i_{j,d,s}}(t-k_s),
\end{align*}
where $J$ is a positive integer, $c_{j,d}(k_1,\ldots,k_d)\in\mathbb R$. This process is called the Volterra series \citep[Equation 2.9.14]{brillinger1981}. Assumption \ref{as_moment} holds true for this process $\bm{X}(t)$ since the mixing property of $\bm{X}(t)$ is inherited from $\bm{\upsilon}(t)$ (see \citealt[p.349]{fz10}) and therefore, the conditions of Lemma \ref{summability_mixing} are fulfilled.
}\end{ex}

Under Assumption \ref{as_moment}, the spectral density matrix $\bm{f}(\lambda):= \skakko{f_{ij}(\lambda)}_{i,j=0,\ldots,K}$ for the  process $(X_0(t),X_1(t),\ldots,X_K(t))$ exists and twice continuously differentiable. For $i\in\{0,\ldots,K\}$, $f_{ii}(\lambda)$ is an auto-spectrum for the process $X_i(t)$ and,  for $i_1,i_2(\neq i_2)\in\{0,\ldots,K\}$, $f_{i_1i_2}(\lambda)$ is a cross-spectrum for the processes $X_{i_1}(t)$ and $X_{i_2}(t)$.

The next theorem shows the model \eqref{model_K} has the following orthogonal representation.

\begin{theorem}\label{thm_decom_in_the_case_of_K}
Suppose that $\bm{f}_K(\lambda):=\skakko{f_{ij}(\lambda)}_{i,j=1,\ldots,K}$ is non-singular for all $\lambda\in[-\pi,\pi]$ and 
${
{\rm det}\skakko{\bm{f}_{j-1}(\lambda)}f_{j0}(\lambda)
\neq
\sum_{i=1}^{j-1}{\rm det}\skakko{\overline{\bm{f}_{i,j}^\ddag(\lambda)}}f_{i0}(\lambda)}$ for $j\in\{2,\ldots,K-1\}$ and for all $\lambda\in[-\pi,\pi]$, where 
$$
{\bm{f}}_{i,j}^\ddag(\lambda):=\skakko{{\bm{f}_{j-1,1}^\flat(\lambda)},\ldots,{\bm{f}_{j-1,i-1}^\flat(\lambda)},{\bm{f}_{j-1,j}^\flat(\lambda)},{\bm{f}_{j-1,i+1}^\flat(\lambda)},\ldots,{\bm{f}_{j-1,j-1}^\flat(\lambda)}}
$$
with $\bm{f}_{a,b}^\flat:=({f}_{1b}(\lambda),\ldots,{f}_{ab}(\lambda))^\top$. There uniquely exists the processes $\{G_j(t)\}$ for $j=1,\ldots,K$, which taking the form of
\begin{align}\label{Gj}
G_j(t):=\sum_{d=1}^j\sum_{k_d=-\infty}^\infty a_{jd}(k_d)X_d(t-k_d)
\end{align}
such that {$X_0(t)=\zeta + \sum_{j=1}^KG_j(t)+\epsilon(t)$}, $G_i$ and $G_j$ are orthogonal for any $i,j(\neq i)\in\{1,\ldots,K\}$, i.e., ${\rm E}G_i(t)G_j(t^\prime)=0$ for any $t,t^\prime\in\mathbb Z$, $\sum_{k=-\infty}^\infty|a_{jd}(k)|<\infty$ for any $j\in \{1,\ldots,K\}$ and $d\in\{1,\ldots,j\}$, and the transfer function $A_{jj}\skakko{e^{-\mathrm{i}\lambda}}:=\sum_{k=-\infty}^\infty a_{jj}(k)e^{-\mathrm{i}k\lambda}$ satisfies $A_{jj}\skakko{e^{-\mathrm{i}\lambda}}\neq0$ for all $\lambda\in[-\pi,\pi]$ and $j\in\{1,\ldots,K-1\}$.
The unique expressions of $G_j$ is given by \eqref{Gj}, where $ a_{jd}$ is determined through the transfer functions $A_{jd}(e^{-\mathrm{i}\lambda})$, as
\begin{align}\label{A11}
A_{11} \skakko{e^{i\lambda}}
:=\frac{f_{10}(\lambda)}{f_{11}(\lambda)}\quad \text{for $j=d=1$},
\end{align}
for $j=d\in\{2,\ldots,K\}$, 
\begin{align}
\label{Ajj}
A_{jj}  \skakko{e^{i\lambda}}
:=
\frac{-\sum_{i=1}^{j-1}{\rm det}\skakko{\overline{\bm{f}_{i,j}^\ddag(\lambda)}}f_{i0}(\lambda)+{\rm det}\skakko{\bm{f}_{j-1}(\lambda)}f_{j0}(\lambda)}{{\rm det}\skakko{\bm{f}_{j}(\lambda)}},
\end{align}
and, for $j\in\{2,\ldots,K\}$ and $d\in\{1,\ldots,j-1\}$,
\begin{align}\label{Ajd}
A_{jd} \skakko{e^{i\lambda}}
:=&-\frac{{\rm det}\skakko{\bm{f}_{d,j}^\ddag(\lambda)}}{{\rm det}\skakko{\bm{f}_{j-1}(\lambda)}}A_{jj} \skakko{e^{i\lambda}}.
\end{align}

\end{theorem}

From Theorem \ref{thm_decom_in_the_case_of_K}, the spectral density for $G_j$ is given by
\begin{align*}
f_{G_1G_1}(\lambda)=
\frac{\left|f_{10}(\lambda)\right|^2}{f_{11}(\lambda)}\quad \text{and }
f_{G_jG_j}(\lambda)=
\frac{\left|A_{jj}  \skakko{e^{i\lambda}}\right|^2{\rm det}\skakko{\bm{f}_{j}(\lambda)}}{{\rm det}\skakko{\bm{f}_{j-1}(\lambda)}} \quad\text{for $j=2,\ldots,K$.}
\end{align*}
By construction, $G_j$ can be interpreted as a process after eliminating the effect of the processes $G_1,\ldots,G_{j-1}$. Thus, we refer to $f_{G_jG_j}$ as a {{\it residual spectral density of order $j$}}.  {These spectra are related to the important indexes.}

{When $K=1$, the squared coherence $\mathcal{C}_1(\lambda)$ is defined, for $\lambda\in[-\pi,\pi]$, as
	\begin{align*}
		\mathcal{C}_1(\lambda):=\frac{f_{G_1G_1}(\lambda)}{f_{00}(\lambda)},
\end{align*}
which is the time series extension of the correlation, measures linear dependence between $\{X_0(t)\}$ and $\{X_1(t)\}$ (see, e.g., \citealt[p.436]{brockwellanddavis1991}).}

{In the case of $K=2$, $G_2$ has the spectral density of the form 
\begin{align*}
	f_{G_2G_2}(\lambda)=
	\frac{\left|f_{11}(\lambda)f_{02}(\lambda)-f_{12}(\lambda) f_{01}(\lambda)\right|^2}{f_{11}(\lambda)\skakko{{f_{11}}(\lambda) { f_{22}}(\lambda)-\left|f_{12}(\lambda)\right|^2}}.
\end{align*}}
{In particular, when $X_{2}(t)$ is defined as the {\it lagged process} $X_2(t):=X_1(t)X_1(t-u)$, the lagged spectral coherence $\mathcal{C}_2(\lambda,u)$ is defined, for $\lambda\in[-\pi,\pi]$, as
\begin{align*}
\mathcal{C}_2(\lambda,u):=\mathcal{C}_1(\lambda)
+\frac{f_{G_2G_2}(\lambda,u)}{f_{00}(\lambda)}.
\end{align*}}
{The lagged coherence $\mathcal{C}_2(\lambda,u)$ measures the nonlinear (quadratic) dependence between $\{X_0(t)\}$ and $\{X_1(t)\}$. See details in \cite{kedem75}. 
In this case, we also refer to $f_{G_2G_2}$ as a {\it lagged spectral density}. }

As a higher-order extension of squared coherence, it is natural to define a {\it squared coherence of order $d$} as 
\begin{align*}
\mathcal{C}_d(\lambda):=\mathcal{C}_{d-1}(\lambda)
+\frac{f_{G_dG_d}(\lambda)}{f_{00}(\lambda)}
=\frac{\sum_{i=1}^{d}f_{G_iG_i}(\lambda)}{f_{00}(\lambda)}.
\end{align*}
Since {$f_{00}(\lambda)=\sum_{i=1}^Kf_{G_iG_i}(\lambda)+{{\rm Var}(\epsilon(t))}/{(2\pi)}$}, $\mathcal{C}_d(\lambda)\in(0,1)$ for any $\lambda\in[-\pi,\pi]$.

The residual spectral density is closely related to the partial cross spectrum, which was introduced by \cite{dahlhaus00}, and is the cross spectrum after removing the linear effect of a certain process (see also \citet[Section 2]{fo11} and \citet[p.969]{eichler08}).
The partial cross spectrum of $X_0$ and $X_2$ given $X_1$ is defined by
\begin{align*}
f_{02|1}(\lambda):=f_{02}(\lambda)-f_{01}(\lambda)f_{11}^{-1}(\lambda)f_{12}(\lambda)
\end{align*}
and we can see the relationship between $f_{G_2G_2}$ and $f_{02|1}$:
\begin{align*}
f_{G_2G_2}(\lambda):=\frac{f_{11}(\lambda)|f_{02|1}(\lambda)|^2}{{f_{11}}(\lambda) { f_{22}}(\lambda)-\left|f_{12}(\lambda)\right|^2}.
\end{align*}
Thus, we can interpret $f_{G_2G_2}$ as the rescaled squared partial cross spectrum of $X_0$ and $X_2$ given $X_1$.

{\begin{rem}
The coherence and the lagged (spectral) coherence can be derived from the model
$$X_0(t):=\sum_{k=-\infty}^\infty b(k)X_1(t-k)+\epsilon(t)$$
and 
\begin{align*}
X_0(t):=&
\sum_{k_1=-\infty}^\infty b_1(k_1)X_1(t-k_1)\\
&+
\sum_{k_2=-\infty}^\infty b_2(k_2)\skakko{X_1(t-k_2)X_1(t-k_2-u)-{\rm E}X_1(t-k_2)X_1(t-k_2-u)}+
\epsilon(t).
\end{align*}
For more details, see \citet[Example 11.6.4]{brockwellanddavis1991} and \cite{kedem75}. In this  sense, our residual spectrum serves as a natural extension of those spectra.
\end{rem}}

\section{Test for the existence of residual spectral density}\label{sec4}

We are interested in the existence of $X_K$ in the model \eqref{model_K}, that is, $b_K(k)=a_{KK}(k)=0$ for any $k\in\mathbb Z$ under the presence of $X_1,\ldots,X_{K-1}$. To this end, we consider the following hypothesis testing problem: the null hypothesis is 
\begin{align}\label{null}
H_0: f_{G_KG_K}(\lambda)=0 \quad\text{$\lambda$-a.e. on $[-\pi,\pi]$}
\end{align}
and the alternative hypothesis is $ K_0: H_0$ does not hold.
The null hypothesis $H_0$ holds true if and only if that $a_{KK}(k)=0$ for all $k\in\mathbb Z$ since 
$a_{KK}(k):=\int_{-\pi}^\pi A_{KK}\skakko{e^{-\mathrm{i}\lambda }}e^{-\mathrm{i}k\lambda}{\rm d}\lambda/(2\pi)$ and 
the proof of Theorem \ref{thm_decom_in_the_case_of_K} shows 
$A_{KK}\skakko{e^{-\mathrm{i}\lambda }}=0$ if and only if 
$
f_{G_KG_K}(\lambda)=0.
$

Suppose the observed stretch $\{X_j(t);t=1,\ldots,n ; j=0,\ldots, K\}$ is available. 
To estimate $\bm{f}(\lambda)$, we define the kernel spectral density estimator $\hat {\bm{f}}(\lambda):=\skakko{\hat{f}_{ij}(\lambda)}_{i,j=0,1,\ldots,K}$ and
\begin{align*}
\hat{f}_{ij}(\lambda):=\frac{1}{2\pi}\sum_{h =1+u-n}^{n-1-u}\omega\skakko{\frac{h}{M_n}}\hat{{\gamma}}_{ij}(h)e^{-\mathrm{i}h\lambda}
\quad
\lambda\in[-\pi,\pi],
\end{align*}
where $M_n$ is a bandwidth parameter satisfies Assumption \ref{as} (A1), $\omega(\cdot)$ is the lag window function defined as $\omega(x):=\int_{-\infty}^\infty W(t) e^{\mathrm{i}xt}\textrm{ d}t$,  $W(\cdot)$ is the kernel function satisfy Assumption \ref{as} (A2), for $h\in\{0,\ldots,n-1-u\}$,
$$\hat{{\gamma}}_{ij}(h):=\frac{1}{n-u-h}\sum_{t=1}^{n-u-h}(X_i{(t+h)}- X_i(.))({ X}_j({t})-{{ X}_j({.})}),$$
and, for $h\in\{-n+1+u,\ldots,-1\}$,
$$\hat{{\gamma}}_{ij}(h):=\frac{1}{n-u+{h}}\sum_{t=-h+1}^{n-u}(X_i{(t+h)}- X_i(.))({ X}_j({t})-{{ X}_j({.})}),$$
with
 ${{X}_i{(.)}}=\sum_{t=1}^{n-u}{X}_i({t})/{{(n-u)}}$  for $i=0,1,\ldots,K$.

\begin{assumption}\label{as}
\begin{enumerate}
\item[(A1)]
A bandwidth parameter satisfies $n/M_n^{9/2}\to0$ and $n/M_n^3\to\infty$ as $n\to\infty$.

\item[(A2)]
The kernel function $W(\cdot)$ is a bounded, nonnegative, even, real, Lipschitz continuous  function such that $\int_{-\infty}^\infty W(t)\textrm{ d}t=1$, $\int_{-\infty}^\infty t^2W(t)\textrm{d}t<\infty$, $\limsup_{t\to\infty}t^2W(t)=0$, $\eta_{\omega,2}:=\int_{-\infty}^\infty \omega^2(x)\textrm{d}x<\infty$, $\eta_{\omega,4}:=\int_{-\infty}^\infty\omega^4(x)\textrm{d}x<\infty$, and
$\sum_{h=-n+1}^{n-1}\omega\skakko{{h}/{M_n}}=O(M_n)$.
\end{enumerate}
\end{assumption}

The condition (A1) is slight stronger than \citet[Assumption 3.3 (iii)]{eichler08}, that is, $n/M_n^{9/2}\to0$ and 
$n/M_n^2\to\infty$ as $n\to\infty$. This is because that he dealt with the zero mean process. On the other hand, we need to estimate ${\rm E}X_i^\prime$ in the kernel density estimator. In the same reason, we assume the condition $\sum_{h=-n+1}^{n-1}\omega\skakko{{h}/{M_n}}=O(M_n)$ in (A2).

From the fact $A_{KK}\skakko{e^{-\mathrm{i}\lambda }}=0$ if and only if $f_{G_KG_K}(\lambda)=0$, {together with the non-singularity of $\bm{f}_K(\lambda)$ and (\ref{Ajj}), the hypothesis (\ref{null}) is equivalent to 
\begin{align*}
	H_0: \Phi_K\skakko{ {\bm{f}}(\lambda)}=0 \quad\text{$\lambda$-a.e. on $[-\pi,\pi]$},
\end{align*}
and subsequently equivalent to 
\begin{align*}
	H_0: \bigintssss_{-\pi}^\pi\left|\Phi_K\skakko{ {\bm{f}}(\lambda)}\right|^2=0,
\end{align*}
where
\begin{align*}
	\Phi_K\skakko{\bm{f}(\lambda)}:=&-\sum_{i=1}^{K-1}{\rm det}\skakko{\overline{ \bm{f}_{i,K}^\ddag(\lambda)}} f_{i0}(\lambda)+{\rm det}\skakko{\bm{f}_{K-1}(\lambda)} f_{K0}(\lambda).
\end{align*}}
Our proposed test statistic is defined as
\begin{align}\label{teststat}
T_n:=\frac{n}{\sqrt {M_n}}\bigintssss_{-\pi}^\pi \left|\Phi_K\skakko{\hat {\bm{f}}(\lambda)}\right|^2{\rm d}\lambda
-  \hat \mu_{n,K},
\end{align}
where 
\begin{align*}
\hat \mu_{n,K}:=&{\sqrt M_n}
\eta_{\omega,2}\bigints_{-\pi}^{\pi}{\rm tr}\skakko{\left.\frac{\partial\Phi_K\skakko{\hat {\bm{f}}(\lambda)}}{\partial\bm Z^\top}\right|_{{\bm Z}=\hat{\bm{f}}(\lambda)}\hat{\bm{f}}(\lambda)\overline{\left.\frac{\partial\Phi_K\skakko{\hat {\bm{f}}(\lambda)}}{\partial\bm Z}\right|_{{\bm Z}=\hat{\bm{f}}(\lambda)}}\hat{\bm{f}}(\lambda)}{\rm d}\lambda,\\
\hat{\bm{f}}_{i,j}^\ddag(\lambda):=&\skakko{{\hat{\bm{f}}_{j-1,1}^\flat(\lambda)},\ldots,{\hat{\bm{f}}_{j-1,i-1}^\flat(\lambda)},{\hat{\bm{f}}_{j-1,j}^\flat(\lambda)},{\hat{\bm{f}}_{j-1,i+1}^\flat(\lambda)},\ldots,{\hat{\bm{f}}_{j-1,j-1}^\flat(\lambda)}},\\
\hat{\bm{f}}_{a,b}^\flat:=&(\hat{f}_{1b}(\lambda),\ldots,\hat{f}_{ab}(\lambda))^\top, \quad \text{and }
\hat{\bm{f}}_{K-1}(\lambda):=\skakko{\hat{f}_{ij}(\lambda)}_{i,j=1,\ldots,K-1}.
\end{align*}
The term $\Phi_K\skakko{\hat {\bm{f}}(\lambda)}$ and $\hat \mu_{n,K}$ are the {estimated} numerator of $A_{KK}$ and a bias correction, respectively. The next theorem shows the asymptotic null distribution of $T_n$.
\begin{theorem}\label{Tn_dist}
Suppose Assumptions \ref{as_moment} and \ref{as}, $\bm{f}_K(\lambda):=\skakko{f_{ij}(\lambda)}_{i,j=1,\ldots,K}$ is non-singular for all $\lambda\in[-\pi,\pi]$ and 
${
{\rm det}\skakko{\bm{f}_{j-1}(\lambda)}f_{j0}(\lambda)
\neq
\sum_{i=1}^{j-1}{\rm det}\skakko{\overline{\bm{f}_{i,j}^\ddag(\lambda)}}f_{i0}(\lambda)}$ for $j\in\{2,\ldots,K-1\}$ and for all $\lambda\in[-\pi,\pi]$. It holds that, under the null $H_0$, $T_n$ defined in \eqref{teststat} converges in distribution to the centered normal distribution with variance $\sigma_K^2$ as $n\to\infty$, where 
\begin{align*}
\sigma_K^2:=4\pi\eta_{\omega,4}
\bigints_{-\pi}^{\pi}&
\left|{\rm tr}\skakko{\left.\frac{\partial\Phi_K\skakko{ {\bm{f}}(\lambda)}}{\partial\bm Z^\top}\right|_{{\bm Z}={\bm{f}}(\lambda)}{\bm{f}}(\lambda)\overline{\left.\frac{\partial\Phi_K\skakko{ {\bm{f}}(\lambda)}}{\partial\bm Z}\right|_{{\bm Z}={\bm{f}}(\lambda)}}{\bm{f}}(\lambda)}\right|^2\\
&+
\left|{\rm tr}\skakko{\left.\frac{\partial\Phi_K\skakko{ {\bm{f}}(\lambda)}}{\partial\bm Z^\top}\right|_{{\bm Z}={\bm{f}}(\lambda)}{\bm{f}}(\lambda)\left.\frac{\partial\Phi_K\skakko{{\bm{f}}(\lambda)}}{\partial\bm Z^\top}\right|_{{\bm Z}={\bm{f}}(\lambda)}{\bm{f}}(\lambda)}\right|^2{\rm d}\lambda.
\end{align*}

\end{theorem}

\begin{rem}{\rm
{For $K=1,2,$ $\Phi_K\skakko{\bm{f}(\lambda)}$ is given by
\begin{align*}
	\Phi_1\skakko{\bm{f}(\lambda)}:=&f_{10}(\lambda)\quad\text{and }\quad\Phi_2\skakko{\bm{f}(\lambda)}:=f_{21}(\lambda)f_{10}(\lambda)-f_{11}(\lambda)f_{20}(\lambda),
\end{align*}
respectively.
Additionally,} the concrete forms of the bias and variance terms are given by, for $K=1$, 
\begin{align*}
\mu_1:=&{\sqrt M_n}
\eta_{\omega,2}\int_{-\pi}^{\pi}f_{11}(\lambda)f_{00}(\lambda){\rm d}\lambda\quad\text{and }\sigma_1^2:=4\pi\eta_{\omega,4}\int_{-\pi}^{\pi}\left| f_{11}(\lambda)f_{00}(\lambda)\right|^{2}{\rm d}\lambda
\end{align*}
and, {for $K=2$,
\begin{align*}
	\mu_2:=&{\sqrt M_n}
	\eta_{\omega,2}\int_{-\pi}^{\pi}\left(f_{11}(\lambda)f_{00}(\lambda)-\left|f_{10}(\lambda)\right|^2\right)\left(f_{11}(\lambda)f_{22}(\lambda)-\left|f_{12}(\lambda)\right|^2\right){\rm d}\lambda
\end{align*}
and
\begin{align*}	
	 \sigma_2^2:=4\pi\eta_{\omega,4}\int_{-\pi}^{\pi}\left|\left(f_{11}(\lambda)f_{00}(\lambda)-\left|f_{10}(\lambda)\right|^2\right)\left(f_{11}(\lambda)f_{22}(\lambda)-\left|f_{12}(\lambda)\right|^2\right)\right|^2{\rm d}\lambda.
\end{align*}}
For $K\geq 3$,
\begin{align*}
\mu_K:=&{\sqrt M_n}
\eta_{\omega,2}\int_{-\pi}^{\pi}(\det{\bm{f}_{K-1}(\lambda)})^{2}
\skakko{f_{00}(\lambda)-\overline{\hat{\bm{f}}_{K-1,0}^{\flat\top}}(\lambda)\bm{f}_{K-1}^{-1}(\lambda){\hat{\bm{f}}_{K-1,0}^\flat}(\lambda)}\\
&\quad\quad\quad\quad\quad\times
\skakko{f_{KK}(\lambda)-\overline{\hat{\bm{f}}_{K-1,K}^{\flat\top}}(\lambda)
\bm{f}_{K-1}^{-1}(\lambda)\hat{\bm{f}}_{K-1,K}^\flat(\lambda)}{\rm d}\lambda\\
\text{and }
\sigma_K^2:=&4\pi\eta_{\omega,4}\int_{-\pi}^{\pi}\left| (\det{\bm{f}_{K-1}(\lambda)})^{2}
\skakko{f_{00}(\lambda)-\overline{\hat{\bm{f}}_{K-1,0}^{\flat\top}}(\lambda)\bm{f}_{K-1}^{-1}(\lambda){\hat{\bm{f}}_{K-1,0}^\flat}(\lambda)}\right.\\
&\quad\quad\quad\quad\quad\times\left.
\skakko{f_{KK}(\lambda)-\overline{\hat{\bm{f}}_{K-1,K}^{\flat\top}}(\lambda)
\bm{f}_{K-1}^{-1}(\lambda)\hat{\bm{f}}_{K-1,K}^\flat(\lambda)
}\right|^{2}{\rm d}\lambda.
\end{align*}}
\end{rem}

A consistent estimator  $\hat\sigma^2_K$ of $\sigma^2_K$ can be constructed by
\begin{align*}
\hat\sigma_K^2:=4\pi\eta_{\omega,4}
\bigints_{-\pi}^{\pi}&
\left|{\rm tr}\skakko{\left.\frac{\partial\Phi_K\skakko{\hat {\bm{f}}(\lambda)}}{\partial\bm Z^\top}\right|_{{\bm Z}=\hat{\bm{f}}(\lambda)}\hat{\bm{f}}(\lambda)\overline{\left.\frac{\partial\Phi_K\skakko{\hat {\bm{f}}(\lambda)}}{\partial\bm Z}\right|_{{\bm Z}=\hat{\bm{f}}(\lambda)}}\hat{\bm{f}}(\lambda)}\right|^2\\
&+
\left|{\rm tr}\skakko{\left.\frac{\partial\Phi_K\skakko{\hat {\bm{f}}(\lambda)}}{\partial\bm Z^\top}\right|_{{\bm Z}=\hat{\bm{f}}(\lambda)}\hat{\bm{f}}(\lambda)\left.\frac{\partial\Phi_K\skakko{\hat {\bm{f}}(\lambda)}}{\partial\bm Z^\top}\right|_{{\bm Z}=\hat{\bm{f}}(\lambda)}\hat{\bm{f}}(\lambda)}\right|^2{\rm d}\lambda
\end{align*}
since $\max_{\lambda\in[-\pi,\pi]}\left\|\hat {\bm{f}}(\lambda)-\bm{f}(\lambda)\right\|_2$, where $\|\cdot\|_2$ denotes the Euclid norm,  converges in probability to zero as $n\to\infty$  (see, e.g., \citealp[Theorem 2.1]{robinson91}). 
From Theorem \ref{Tn_dist}, the test which rejects $H_0$ whenever $T_n/{\hat\sigma_K}\geq z_{\alpha}$, where $z_{\alpha}$ denotes the upper $\alpha$-percentile of the standard normal distribution, has asymptotically size $\alpha$. The following theorem shows the power of the test tends to one as $n\to\infty$.

\begin{theorem}\label{thm_cons}
Suppose Assumptions \ref{as_moment} and \ref{as}, $\bm{f}_K(\lambda):=\skakko{f_{ij}(\lambda)}_{i,j=1,\ldots,K}$ is non-singular for all $\lambda\in[-\pi,\pi]$ and 
${
{\rm det}\skakko{\bm{f}_{j-1}(\lambda)}f_{j0}(\lambda)
\neq
\sum_{i=1}^{j-1}{\rm det}\skakko{\overline{\bm{f}_{i,j}^\ddag(\lambda)}}f_{i0}(\lambda)}$ for $j\in\{2,\ldots,K-1\}$ and for all $\lambda\in[-\pi,\pi]$.  The test based on $T_n$ is consistent.
\end{theorem}

For the case that $K=2$ and we choose the lagged process $X_2(t):=X_1(t)X_1(t-u)$ as a second covariate process, the lag $u$ must be determined by a data analyst.  In order to to choose the 
plausible lag, we propose the following criterion: denote the spectral density matrix $\bm{f}(\lambda)$ as $\bm{f}(\lambda,u)$ to emphasize the dependence on $u$. For some $L\in\mathbb N$, we define
\begin{align*}
\hat u
:= \argmax_{u\in\{0,\ldots,L\}}\int_{-\pi}^\pi q\skakko{\hat {\bm{f}}(\lambda,u)}{\rm d}\lambda,
\end{align*}
where $q$ is some function. The intuitive choice of $q$ would be $q\skakko{\hat {\bm{f}}(\lambda)}:=\hat f_{G_2G_2}(\lambda,u)$, where 
\begin{align*}
\hat f_{G_2G_2}(\lambda,u):=
\frac{\left|\hat f_{11}(\lambda)\hat f_{20}(\lambda,u)-\hat f_{21}(\lambda,u)\hat  f_{10}(\lambda)\right|^2}{\hat f_{11}(\lambda)\skakko{{\hat f_{11}}(\lambda) { \hat f_{22}}(\lambda,u)-\left|\hat f_{21}(\lambda,u)\right|^2}}.
\end{align*}
 This choice is minimizing the mean square error ${\rm E}\epsilon_t^2$ since ${\rm E}\epsilon_t^2=\int_{-\pi}^\pi f_{ZZ}(\lambda)-f_{G_1G_1}(\lambda)-f_{G_2G_2}(\lambda,u) {\rm d}\lambda$. 
If a data analyst is interested in the maximum peak of $f_{G_2G_2}(\lambda,u)/f_{00}(\lambda)$, we can choose $q$ as
\begin{align*}
q\skakko{\hat {\bm{f}}(\lambda)}
:=
\delta\skakko{\lambda-\argmax_{\omega\in[-\pi,\pi]}\skakko{\frac{\hat f_{G_2G_2}(\omega,u)}{\hat f_{00}(\omega)}}}
\frac{\hat f_{G_2G_2}(\omega,u)}{\hat f_{00}(\omega)},\end{align*}
where $\delta$ is the Dirac delta function,  
which corresponds to the choice based on the residual coherence \citep{kkk14}. 
Similarly, for the case that $K$ is greater than 2 and we choose the covariate $X_K(t)$ depends on the lag $u_1,\ldots,u_{K-1}$, e.g., $\prod_{j=1}^{K-1} X_j(t-u_j)$, we propose the criterion 
\begin{align*}
(\hat u_1,\ldots, \hat u_{K-1})
:= \argmax_{(u_1,\ldots,u_{K-1})\in\{0,\ldots,L\}^{K-1}}\int_{-\pi}^\pi q\skakko{\hat {\bm{f}}(\lambda,u_1,\ldots,u_{K-1})}{\rm d}\lambda.
\end{align*}

\begin{rem}\label{rem3.2}{\rm 
One may be skeptical of handling $\hat u_n$ as a known constant in the testing procedure. In that case, we can consider the following hypothesis, for candidates of lags $u_1,\ldots,u_L\in\mathbb N\cup \{0\}$,
\begin{align}\label{mult_null}
\tilde H_0: f_{G_KG_K}(\lambda,u)=0 \quad\text{$\lambda$-a.e. on $[-\pi,\pi]$ for all $u\in\{u_1,\ldots,u_L\}$.}
\end{align}
and $\tilde K_0: \tilde H_0$ does not hold, where $L$ is a given constant, and the corresponding test statistic is defined by
\begin{align}\label{mult_teststat}
\tilde T_n:=\frac{n}{\sqrt {M_n}}\int_{-\pi}^\pi \left\|\bm {\Phi}\skakko{\tilde{\bm{f}}(\lambda)}\right\|_2^2{\rm d}\lambda
-\hat {\tilde\mu}_{n,K},
\end{align}
where $\|\cdot\|_2$ denotes the Euclid norm,  
$\tilde{\bm{f}}(\lambda)$ is the spectral density matrix of

\noindent
$(X_0(t),X_1(t),\ldots,X_{K-1}(t),X_{K,u_1}(t),\ldots,X_{K,u_L}(t))$, $X_{K,u_j}(t)$ is the $K$-th covariate corresponding to the lag $u_j$, 
$$
\bm {\Phi}\skakko{\tilde{\bm{f}}(\lambda)}:=\skakko{
\Phi_{K,u_1}\skakko{\tilde{\bm{f}}(\lambda)}
,\ldots,
\Phi_{K,u_L}\skakko{\tilde{\bm{f}}(\lambda)}
}^\top,\quad
\Phi_{K,u_j}\skakko{\tilde{\bm{f}}(\lambda)}
:=\Phi_{K}\skakko{\hat {\bm{f}}(\lambda,u_j)},
$$
$$
\hat {\tilde\mu}_{n,K}:=
\sqrt {M_n}\eta_{\omega,2}
\int_{-\pi}^\pi{\rm tr}
\mkakko{
\Gamma_{\bm {\Phi}}(\lambda)
\skakko{\tilde{\bm{f}}^\top(\lambda)\otimes \tilde{\bm{f}}(\lambda)}
}{\rm d}\lambda,
$$
$\hat {\bm{f}}(\lambda,u_j)$ is the kernel density estimator for the spectral density matrix of the process $(X_0(t),X_1(t),\ldots,X_{K-1}(t),X_{K,u_j}(t))$, and
\begin{align*}
\Gamma_{\bm {\Phi}}(\lambda)
:= \sum_{k=1}^L {\rm vec}
\left.
\skakko{\frac{\overline{\partial{\Phi_{K,u_k}}\skakko{{\bm Z}}}}{\partial {\bm Z}}}
{\rm vec}\skakko{\frac{{\partial {\Phi_{K,u_k}}\skakko{{\bm Z}}}}{\partial {\bm Z}}}^\top
\right|_{{\bm Z}=\tilde{\bm{f}}(\lambda)}.
\end{align*}
 Then, $\tilde T_n$ converges in distribution to the centered  normal distribution with variance $\tilde\sigma^2_K$ as $n\to\infty$, where
\begin{align*}
\tilde\sigma^2_K:=&
4\pi \eta_{\omega,4}
\int_{-\pi}^\pi{\rm tr}
\mkakko{
\Gamma_{\bm {\Phi}}(\lambda)
\skakko{\tilde{\bm{f}}^\top(\lambda)\otimes \tilde{\bm{f}}(\lambda)}
\skakko{\Gamma_{\bm {\Phi}}(\lambda)
+\Gamma_{\bm {\Phi}}^\top(-\lambda)
}
\skakko{\tilde{\bm{f}}^\top(\lambda)\otimes \tilde{\bm{f}}(\lambda)}
}{\rm d}\lambda.
\end{align*}
Then, the test which rejects $\tilde H_0$ whenever $T_n/{\hat{\tilde{\sigma}}_K}\geq z_{\alpha}$ has asymptotically size $\alpha$ and is consistent, where $\hat{\tilde{\sigma}}_K$ is a consistent estimator of $\tilde\sigma_K$. The proof is omitted since it is analogous to the proof of Theorems \ref{Tn_dist} and \ref{thm_cons}. 
}\end{rem}

One may be interested in the parameters $\{b_i(k);k\in\mathbb Z, i\in\{1,\ldots,K\}\}$. Before closing this section, we construct estimators of $\{b_i(k);k\in\mathbb Z, i\in\{1,\ldots,K\}\}$. The estimator of $\{b_i(k)$ can be defined, for any $k\in\mathbb Z$, by
$
\hat b_i(k):=\sum_{j=i}^K\hat a_{ji}(k),
$ where, for $i,j(\leq i) $, 
$\hat a_{ij}(k):=\int_{-\pi}^\pi \hat A_{ij}\skakko{e^{-\mathrm{i}\lambda }}e^{-\mathrm{i}k\lambda}{\rm d}\lambda/(2\pi)$,
where $\hat A_{ij}\skakko{e^{-\mathrm{i}\lambda }}$ is defined as 
\eqref{A11}--\eqref{Ajd} but the spectral density is replaced with the kernel density estimator.
Then, we have the consistency of estimators.
\begin{theorem}\label{thm_est}
Suppose Assumptions \ref{as_moment} and \ref{as}, $\bm{f}_K(\lambda):=\skakko{f_{ij}(\lambda)}_{i,j=1,\ldots,K}$ is non-singular for all $\lambda\in[-\pi,\pi]$ and 
${
{\rm det}\skakko{\bm{f}_{j-1}(\lambda)}f_{j0}(\lambda)
\neq
\sum_{i=1}^{j-1}{\rm det}\skakko{\overline{\bm{f}_{i,j}^\ddag(\lambda)}}f_{i0}(\lambda)}$ for $j\in\{2,\ldots,K-1\}$ and for all $\lambda\in[-\pi,\pi]$. For any $k\in\mathbb Z$, the estimators $\hat b_i(k)$ for $i \in \{1,\ldots,K\}$ converges in probability to $b_i(k)$ as $n\to\infty$.
\end{theorem}

\section{Numerical study}\label{sec5}

This section presents the finite sample performance of the proposed test. Let $\{\epsilon_j(t);t\in\mathbb Z\}$ be an i.i.d.\ standard normal random variable, where $j=0,\ldots,4$ and each variable is independent of the others, $\{X_j(t);t\in\mathbb Z\}$ for $j=1,\ldots,3$ be the AR(1) model defined as $X_j(t)=0.4X_j(t-1)+\epsilon_j(t)$, and $\{X_4(t);t\in\mathbb Z\}$ be the process given by $X_4(t)=X_2(t)+\epsilon_4(t)$. We consider 14 cases presented in Table \ref{tb1}. Cases 1, 4, and 7 are associated with to the null hypothesis, while the other cases are related to the alternative hypothesis. Cases 1--3, 4--5, and 6--10 are associated with linear terms only, and correspond to tests with $K=1$, $K=2$, and $K=3$, respectively. In contrast, Cases 11--14 involve non-linear (interaction) terms and correspond to tests with $K=2$. It is worth noting that Case 10 is close to the null since our test is designed to remove the effect, in the sense of orthogonality, of the processes that are already included in the model.
\begin{table}[h]
\centering 
\caption{The models and tests for numerical simulations.}
\scalebox{0.8}{
 \begin{tabular}{ccc}
Case   &  model & test\\\hline
1&$X_0(t) = \epsilon_0(t)$&the existence of $X_1$\\
2&$X_0(t) = 0.05 X_1(t) + \epsilon_0(t)$&
the existence of $X_1$\\
3&$X_0(t) = 0.1 X_1(t) + \epsilon_0(t)$&the existence of $X_1$\\\hline
4&$X_0(t) = X_1(t) + \epsilon_0(t)$&the existence of $X_2$ under the presence of $X_1$\\
5&$X_0(t) = X_1(t)+ 0.05 X_2(t) + \epsilon_0(t)$&the existence of $X_2$ under the presence of $X_1$\\
6&$X_0(t) = X_1(t)+ 0.1 X_2(t) + \epsilon_0(t)$&the existence of $X_2$ under the presence of $X_1$\\\hline
7&$X_0(t) = X_1(t) + X_2(t)+ \epsilon_0(t)$&the existence of $X_3$ under the presence of $X_1,X_2$\\
8&$X_0(t) = X_1(t) + X_2(t)+ 0.05 X_3(t) + \epsilon_0(t)$&the existence of $X_3$ under the presence of $X_1,X_2$\\
9&$X_0(t) = X_1(t) + X_2(t)+ 0.1 X_3(t) + \epsilon_0(t)$&the existence of $X_3$ under the presence of $X_1,X_2$\\\hline
10&$X_0(t) = X_1(t) + X_2(t)+ 0.05 X_4(t) + \epsilon_0(t)$&the existence of $X_4$ under the presence of $X_1,X_2$\\\hline
11&$X_0(t) = X_1(t) + 0.05X_1^2(t) + \epsilon_0(t)$&the existence of $X_1(t)^2$ under the presence of $X_1$\\
12&$X_0(t) = X_1(t) + 0.05X_1^2(t) + \epsilon_0(t)$&the existence of $X_1(t)X_1(t-1)$ under the presence of $X_1$\\
13&$X_0(t) = X_1(t) + 0.05X_1^2(t) + \epsilon_0(t)$&the existence of $X_1(t)X_1(t-2)$ under the presence of $X_1$\\
14&$X_0(t) = X_1(t) + 0.05X_1^2(t) + \epsilon_0(t)$&the existence of $X_1(t)X_1(t-3)$ under the presence of $X_1$
\end{tabular}}
\label{tb1}
\end{table}
The sample size and significance level is set to $n=250,500,1000,2000$ and $0.05$, respectively. For each case, we generate a time series with length of $n$ and apply our test. We iterate this procedure 1000 times and calculate empirical size or power of the test.

The results are given in Table \ref{tb2}. The empirical size (corresponds to  Cases 1, 4, 7) is acceptable. The empirical power for Cases 2, 3, 5, 6, 8, 9, 11, increases  as the sample size gets larger. For Case 10, the empirical power is small, as expected. For Cases 12--14, which are the tests for the existence of $X_1(t)X_1(t-1), X_1(t)X_1(t-2), X_1(t)X_1(t-3)$ in the model, respectively, the power is smaller as the lag of the lagged process increases, which is reasonable. Overall, the proposed test shows good performance.

For Case 12, $X_0(t) = X_1(t) + 0.05X_1^2(t) + \epsilon_0(t)$ can be decomposed into three components 

%

\begin{enumerate}
\item[(i)]
{$X_1(t) + 0.05 {\rm P}_{\overline{\rm sp}\{X_1(t);t\in \mathbb Z\}}X_1^2(t)$}
\item[(ii)]
{$0.05{\rm P}_{\overline{\rm sp}\{X_1(t)X_1(t-1);t\in \mathbb Z\}}
\skakko{X_1^2(t)-{\rm P}_{\overline{\rm sp}\{X_1(t);t\in \mathbb Z\}}X_1^2(t)}$}
\item[(iii)]
{$\epsilon_0(t)+  0.05(X_1^2(t)
-{\rm P}_{\overline{\rm sp}\{X_1(t);t\in \mathbb Z\}}X_1^2(t)$
$-{\rm P}_{\overline{\rm sp}\{X_1(t);t\in \mathbb Z\}}(X_1^2(t)-{\rm P}_{\overline{\rm sp}\{X_1(t);t\in \mathbb Z\}}X_1^2(t))$}
\end{enumerate}

\noindent
where ${\rm P}_A$ and $\overline{\rm sp}A$ denote the projection operator onto the set $A$ and the closed span of $A$, respectively. This decomposition tell us that the rejection probability of the test is not $0.05$ when the model includes a second input (in our case, $X_1(t)^2$) that is not the one we wish to test ($X_1(t)X_1(t-1)$) since the second input partially encompasses the input that we intend to test (the term (ii)) unless they are orthogonal, which is not the case for $X_1(t)^2$ and $X_1(t)X_1(t-1)$.
Therefore, Case 12 corresponds to the alternative. Cases 13 and 14 are similar. One may notice that (iii) is not i.i.d.\ disturbance process any more. Since independence can be relaxed to orthogonality between the error term and covariates and it can be time series, our test can be applied.

\begin{table}[h]
\centering 
\caption{Empirical size and power of the proposed test for each case.}
\scalebox{1}{
 \begin{tabular}{ccccc}
Case  $\backslash$ $n$& 250 & 500 & 1000 & 2000\\\hline
1& 0.068& 0.085& 0.066& 0.068\\
2& 0.156& 0.224& 0.333& 0.532\\
3& 0.369& 0.631& 0.865& 0.990\\\hline
4& 0.069& 0.075& 0.069& 0.059\\
5& 0.148& 0.214& 0.284& 0.430\\
6& 0.353& 0.550& 0.768& 0.964\\\hline
7& 0.069& 0.074& 0.061& 0.062\\
8& 0.140& 0.186& 0.291& 0.420\\
9& 0.306& 0.511& 0.736& 0.936\\\hline
10&0.114& 0.13& 0.173& 0.214\\\hline
11& 0.213& 0.277& 0.465& 0.720\\
12& 0.148& 0.212& 0.315& 0.510\\
13& 0.075& 0.108& 0.145& 0.201\\
14& 0.081& 0.075& 0.077& 0.107\\
\end{tabular}}
\label{tb2}
\end{table}

\section{Empirical study}\label{sec6}

\begin{figure}
\begin{center}
\includegraphics[width=\linewidth]{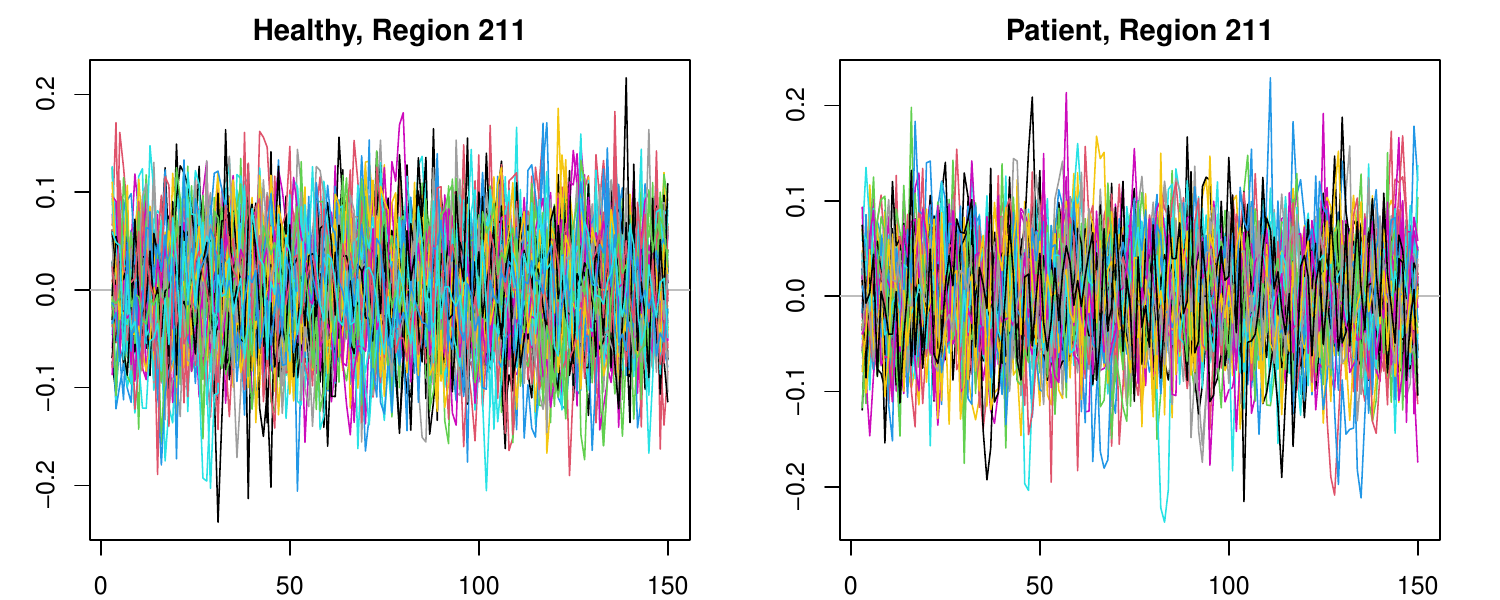}
  \caption{{Plots of time series for 41 healthy individuals (the left panel) and 49 schizophrenia patients (the right panel) for the region 211.}}
  \label{fig:plot_211}
\end{center}
\end{figure}

Our data consist of fMRI time series from 61 healthy individuals (controls) and 49 schizophrenia patients (cases), where in both cases the brain was divided into 246 regions based on the Brainnetome atlas (see \cite{fan16}), each giving a region-level time series of length {148}. The data are described in \cite{culbreth21}. {Figure \ref{fig:plot_211} shows the plots of time series for 61 healthy individuals and 49 schizophrenia patients for the region 211, respectively.}

{The regions 211--246 correspond to subcortical nuclei. Many research has been conducted for the association between schizophrenia and subcortical abnormalities, for example, \cite{fan19}. Thus, it is of particular interest to detect structural differences within these regions between healthy people and patients.}

{We investigate whether time series from region $i\in\{211,\ldots,246\}$ can be explained by time series from the other region $j\in\{211,\ldots,246\}\backslash\{i\}$ through two scenarios. 
In the first scenario (i), we apply our test based on $T_n$ defined in \eqref{teststat} for the hypothesis \eqref{null} with $K=1$ to time series $X_0$ and $X_1$ observed from regions $i\in\{211,\ldots,246\}$ and from $j\in\{211,\ldots,246\}\backslash\{i\}$ for healthy individuals and patients, respectively. 
In the second scenario (ii), we investigate the hypothesis \eqref{null} and apply our test based on $T_n$ with $K=2$ and {$X_{2}(t)=X_{1}(t)X_{1}(t- \hat u)$}, where 
\begin{align*}
\hat u
:= \argmax_{u\in\{0,\ldots,5\}}\int_{-\pi}^\pi \hat f_{G_2G_2}(\lambda,u){\rm d}\lambda
\end{align*}
(see the discussion below Theorem \ref{thm_cons}), to time series $X_0$ and $X_1$ observed from regions $i\in\{163,\ldots,188\}$ and from $j\in\{163,\ldots,188\}\backslash\{i\}$for healthy individuals and patients, respectively. It is assumed that the linear term $\{X_1(t)\}$ has already been incorporated into included in the model.
Note that the scenarios (i) and (ii) correspond to tests for $f_{G_1G_1}(\lambda)=0$ and $f_{G_2G_2}(\lambda)=0$ for $\lambda$ a.e.\ on $[-\pi,\pi]$, respectively.}

{Figure \ref{fig:res_211_246} illustrates the rejection probabilities (regarding individuals) of the test results described in scenario (i) for both healthy individuals and patients, respectively.
The regions can be divided in terms of gyrus regions (See Table \ref{tab:gyrus_regions}). These regions are emphasized by blue lines in Figure \ref{fig:res_211_246}. The regions belonging to the Amygdala and the Hippocampus are well connected, as are the regions belonging to the Thalamus. Additionally, the regions belonging to the Basal Ganglia are well connected, with the exception of regions 227 and 228. The rejection probabilities for regions 227 and 228 in the Basal Ganglia to the other regions are relatively low. Similarly, those for regions 242--244 in the Thalamus to the other regions in the Amygdala, Hippocampus, and Basal Ganglia are relatively low. These regions are emphasized by yellow dotted lines in Figure \ref{fig:res_211_246}.}
\begin{table}[htbp]
\centering
\begin{tabular}{|c|c|}
\hline
\textbf{Gyrus Region} & \textbf{Region Numbers} \\ \hline
\multirow{1}{*}{Amygdala} & 211--214 \\ \hline
\multirow{1}{*}{Hippocampus} & 215--218 \\ \hline
\multirow{1}{*}{Basal Ganglia} & 219--230 \\ \hline
\multirow{1}{*}{Thalamus} & 231--246 \\ \hline
\end{tabular}
\caption{Corresponding Gyrus Regions of Our Dataset}
\label{tab:gyrus_regions}
\end{table}

{Figure \ref{fig:res_211_246_2} presents the corresponding the rejection probability for scenario (ii). The rejection probability for each pixel is explicitly displayed. In each panel, the left and bottom margins display the region numbers $i\in\{211,\ldots,246\}$ corresponding to time series $X_0$ and $X_1$ from brain region $i$, respectively.
Additionally, it appears that the rejection probabilities of inter-gyrus regions for patient tend to be higher than those for healthy individuals.
}

{Table \ref{tab:rejection_probabilities} provides a concise summary comparing the  rejection probabilities between healthy individuals and patients. As expected, Figure \ref{fig:res_211_246} and Table \ref{tab:rejection_probabilities} indicate that the rejection probabilities for patients tend to be lower than those for healthy individuals, consistent with findings in \cite{liu18}, with several pixels yielding a rejection probability of one. Conversely, Figure \ref{fig:res_211_246_2} and Table \ref{tab:rejection_probabilities} demonstrates that the rejection probabilities for patients tend to be higher than those for healthy individuals.
\begin{figure}
\begin{center}
\includegraphics[width=2.7in]{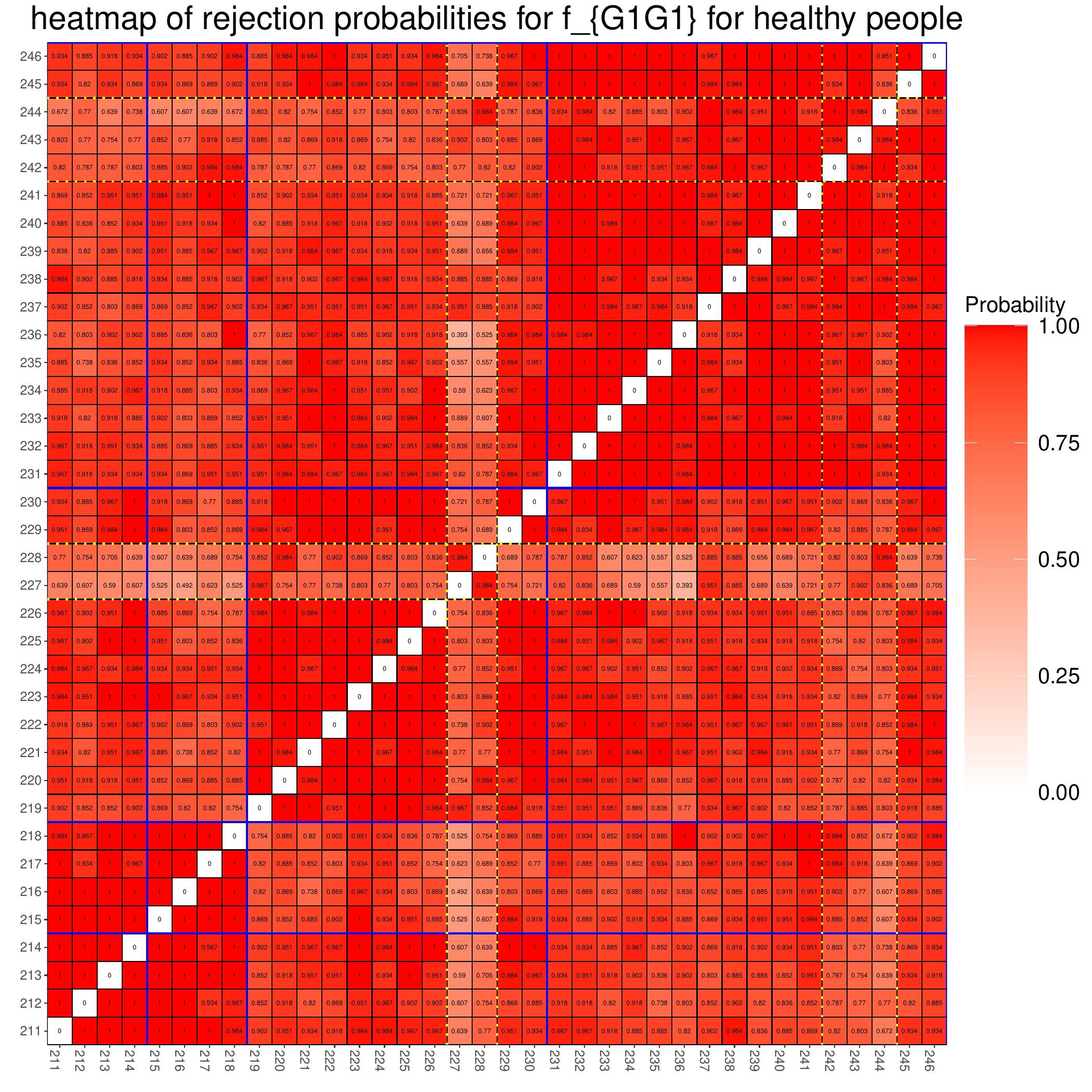}
  \includegraphics[width=2.7in]{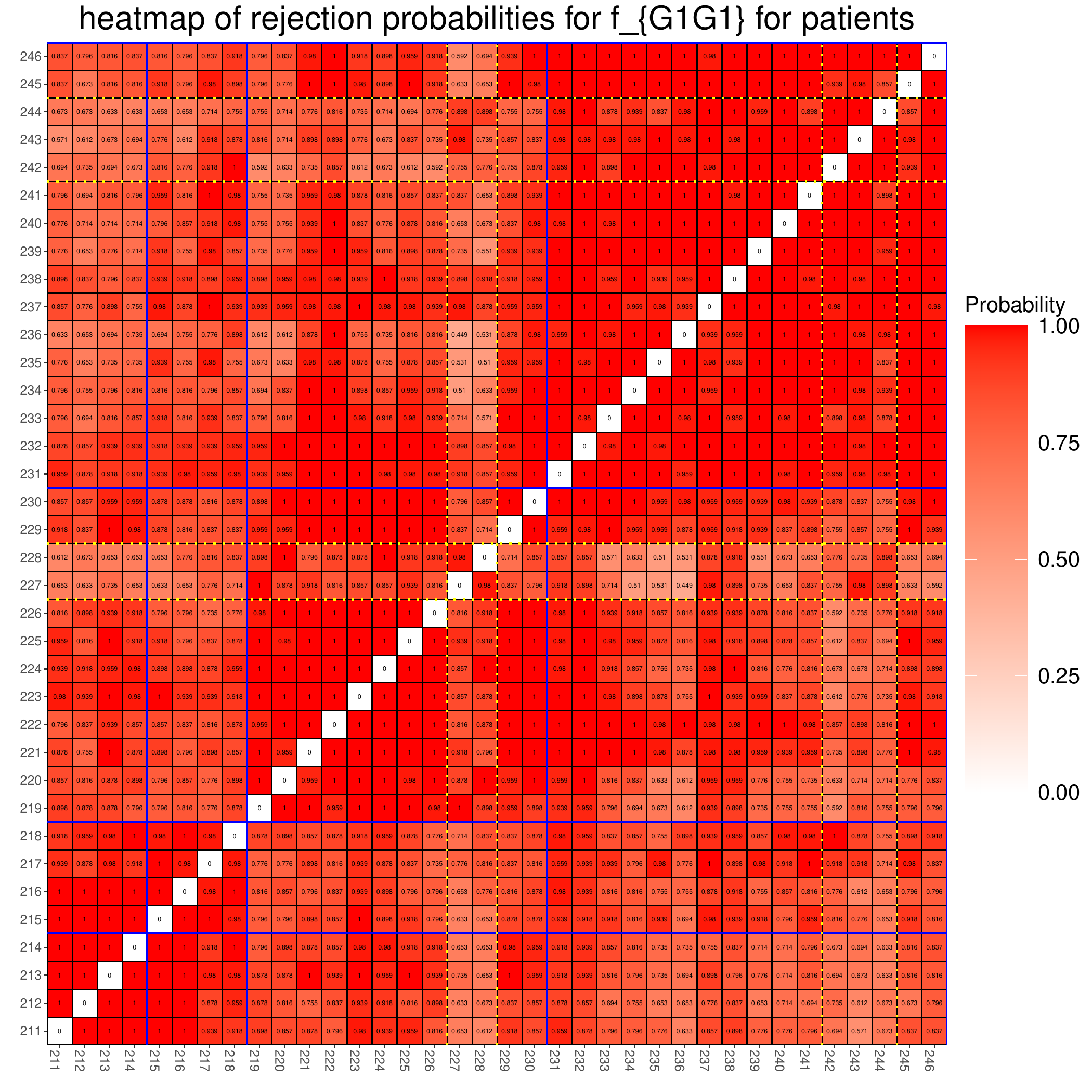}
  \caption{The left and right panels depict the rejection probabilities of the test results outlined in scenario (i) for healthy individuals and patients, respectively. In each panel, the left and bottom margins show the region numbers $i\in\{211,\ldots,246\}$ associated with time series $X_0$ and $X_1$ from the brain region $i$, respectively.}
  \label{fig:res_211_246}
\end{center}
\end{figure}
\begin{figure}
\begin{center}
\includegraphics[width=2.7in]{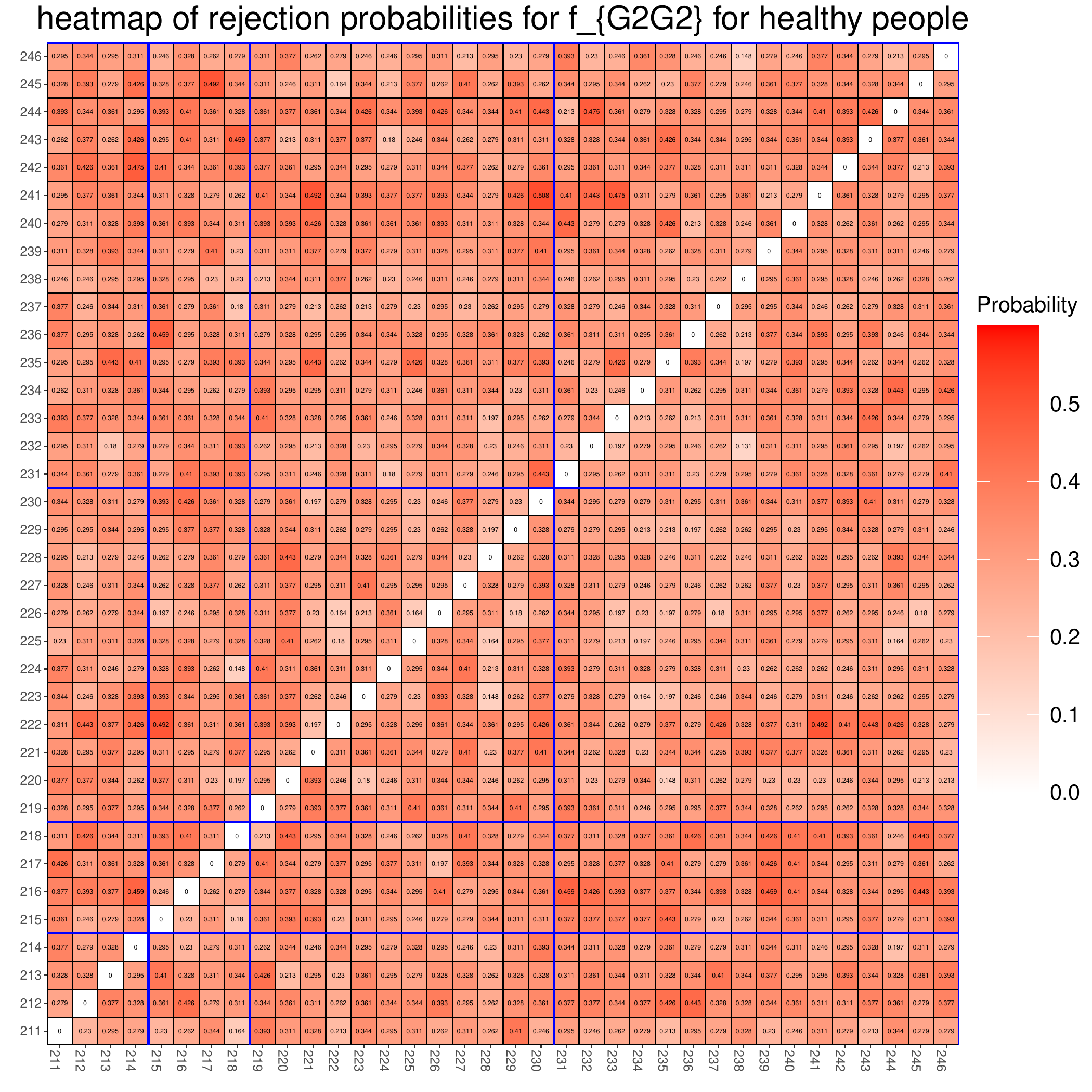} \includegraphics[width=2.7in]{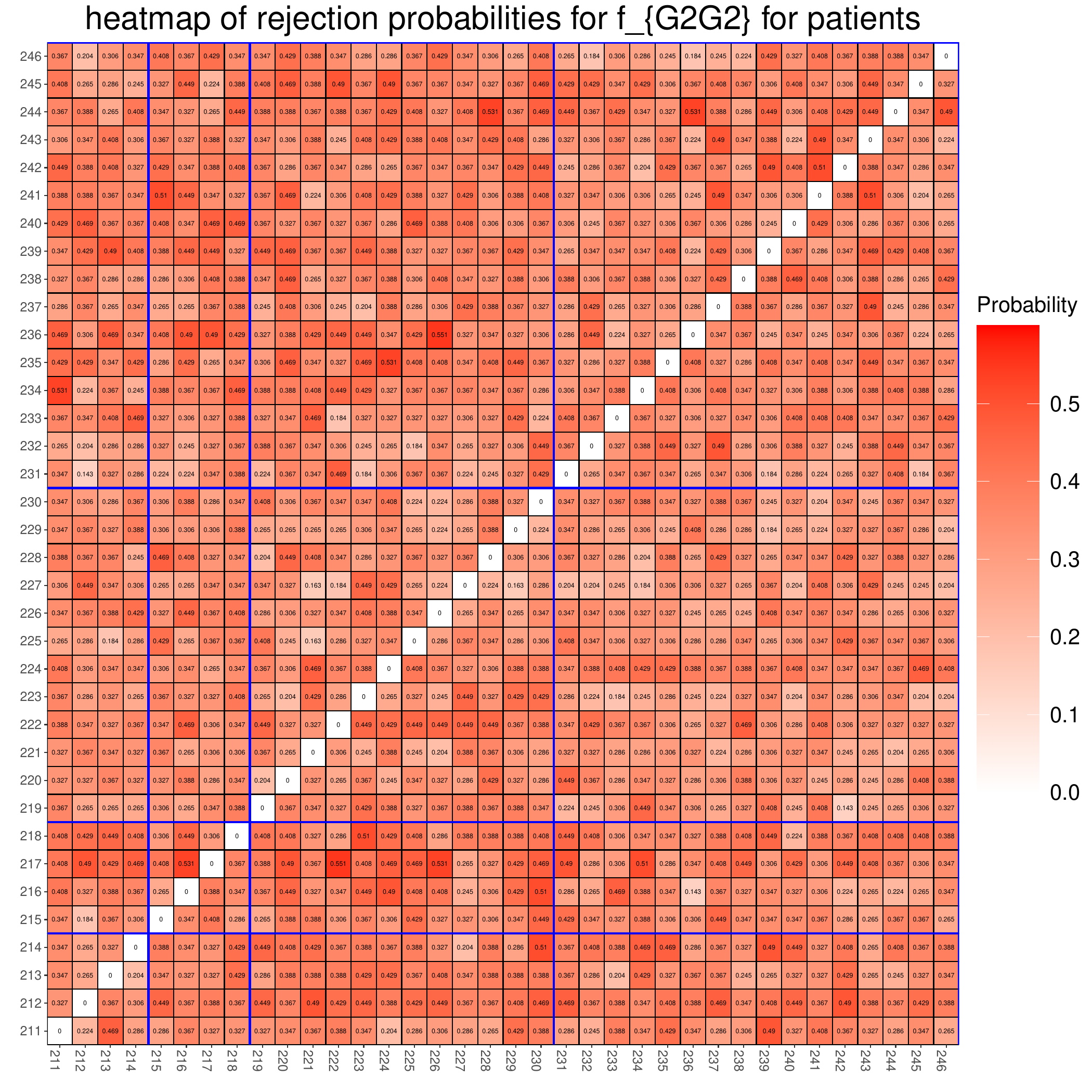}
  \caption{The left and right panels depict the rejection probabilities of the test results outlined in scenario (ii) for healthy individuals and patients, respectively. In each panel, the left and bottom margins show the region numbers $i\in\{211,\ldots,246\}$ associated with time series $X_0$ and $X_1$ from the brain region $i$, respectively.}
  \label{fig:res_211_246_2}
\end{center}
\end{figure}
\begin{table}[htbp]
\centering
\caption{Comparison of Rejection Probabilities between healthy individuals and patients}
\label{tab:rejection_probabilities}
\begin{tabular}{|c|c|c|}
\hline
\multirow{2}{*}{\textbf{Probability}} & \multicolumn{2}{c|}{\textbf{Outcome}} \\ \cline{2-3} 
 & \textbf{a test for $f_{G_1G_1}$} & \textbf{a test for $f_{G_2G_2}$} \\ \hline
\textbf{Patients $>$ Healthy} &  30.00\% & 62.14\% \\ \hline
\textbf{Patients $<$ Healthy} &  48.25\% & 37.86\%\\ \hline
\textbf{Patients $=$ Healthy} &  21.75\% & 0.00\%\\ \hline
\end{tabular}
\end{table}}

{In summary, our findings suggest that there are statistically significant linear relationships between fMRI data from different regions of the brain, consistent with prior research. Furthermore, beyond these linear relationships, our analysis reveals the presence of statistically significant nonlinear associations, which contribute to the explanation of fMRI data beyond the scope of linear relationships alone. Interestingly, our investigation indicates that non-linear functional connectivity  tends to be higher among schizophrenia patients compared to healthy subjects. This contrasts with the findings in the linear case, where patients tend to exhibit lower functional connectivity than healthy subjects when considering only linear measures.}

\begin{appendix}\label{appen}

\section{Proofs}
\subsection{Proof of Lemma 2.1}

First, we derive the upper bound for the quantity $${\rm cum}\{X_{i_1}(0),X_{i_2}(s_2),\ldots, X_{i_\ell}(s_\ell)\}.$$ 
Its derivation is analogous to \citet[Proposition 3.1]{dette16}.
In the case $s_2=\dots=s_\ell=0$, there exists, by the moment assumption, a constant $C_\ell$ such that
\begin{align*}
\sup_{i_1,\ldots,i_\ell\in\{0,1,\ldots,K\}}\abs{{\rm cum}\{X_{i_1}(0),X_{i_2}(0),\ldots, X_{i_\ell}(0)\}}\leq C_\ell.
\end{align*}
We  shall consider the case when there exists $j\in\{0,1,\ldots,K\}$ such that $s_j\neq0$. Let $(s_{(1)},\ldots,s_{(\ell)})$ be the order statistic of $(0, s_2,\ldots,s_{\ell})$ and $j^*$ be an index corresponding to the maximum lag of the order statistic, i.e., $j^*:=\argmax_{j\in\{0,1,\ldots,K-1\}}(s_{(j+1)}-s_{(j)})$. 
Recall the inequalities that
\begin{enumerate}
\item[(i)] for $x_1,\ldots,x_p,y_1\ldots,y_p\in\mathbb R$ whose absolute values are of bounded by $C_{\rm ineq}$,
\begin{align*}
\abs{\prod_{j=1}^px_j-\prod_{j=1}^py_j}\leq C_{\rm ineq}^{p-1}\sum_{j=1}^p\abs{x_j-y_j}
\end{align*}
(see, e.g., \citealt[Lemma F.5]{pt13suppl}),

\item[(ii)] for any $\delta>1, s_1,s_2\in\mathbb Z$, and the alpha-mixing stationary processes $W_1(s_1)$ and $W_2(s_2)$ with $\lceil \frac{2\delta}{\delta-1} \rceil$-th moment,
\begin{align*}
\abs{{\rm Cov}(W_1(s_1), W_2(s_2))}\leq 8 \alpha^{\frac{1}{\delta}}(|s_2-s_1|)\skakko{{\rm E}W_1^{\frac{2\delta}{\delta-1}}(s_1)}^{\frac{\delta-1}{2\delta}} \skakko{{\rm E}W_2^{\frac{2\delta}{\delta-1}}(s_2)}^{\frac{\delta-1}{2\delta}},
\end{align*}
where $\alpha(\cdot)$ is the alpha-mixing coefficient. See, e.g., \citet[Theorem 3 (i), p.9]{doukhan94}.
\end{enumerate}
From the above inequalities (i) and (ii), it follows, for appropriate indexes $i_{(1)},\ldots,i_{(\ell)}\in\{0,1,\ldots,K\}$ and the sets $\upsilon:=\{j\in{1,\ldots,\ell}: s_{(j)}\leq s_{(j^*)}\}$ and $\upsilon^\mathsf{c}:=\{j\in{1,\ldots,\ell}: s_{(j)}>s_{(j^*)}\}$,
\begin{align*}
&\abs{{\rm cum}\{X_{i_1}(0),X_{i_2}(s_2),\ldots, X_{i_\ell}(s_\ell)\}}\\
=&
\abs{{\rm cum}\{X_{i_{(1)}}(s_{(1)}),X_{i_{(2)}}(s_{(2)}),\ldots, X_{i_{(\ell)}}(s_{(\ell)})\}}\\
&-  \abs{{\rm cum}\{X_{i_{(1)}}(s_{(1)}),X_{i_{(2)}}(s_{(2)}),\ldots,X_{i_{(j^*)}}(s_{(j^*)}), X^*_{i_{(j^*+1)}}(s_{(j^*+1)}),\ldots, X^*_{i_{(\ell)}}(s_{(\ell)})\}}\\
\leq&
\Bigg|\sum_{(\nu_1,\ldots,\nu_p)}(-1)^{p-1}(p-1)!\Bigg\{
\skakko{{\rm E}\prod_{j\in\nu_1}X_{i_{(j)}}(s_{(j)})}\ldots\skakko{{\rm E}\prod_{j\in\nu_p}X_{i_{(j)}}(s_{(j)})}\\
&-\skakko{{\rm E}\prod_{j\in\nu_1 \cap \upsilon}X_{i_{(j)}}(s_{(j)})}
\skakko{{\rm E}\prod_{j\in\nu_1 \cap\upsilon^\mathsf{c}}X_{i_{(j)}}^*(s_{(j)})}
 \ldots
\skakko{{\rm E}\prod_{j\in\nu_p \cap \upsilon}X_{i_{(j)}}(s_{(j)})}\\&\quad\times
\skakko{{\rm E}\prod_{j\in\nu_p \cap\upsilon^\mathsf{c}}X_{i_{(j)}}^*(s_{(j)})}
\Bigg\}\Bigg|\\
\leq&
\sum_{(\nu_1,\ldots,\nu_p)}(p-1)!
C_{\rm ineq,\ell}^{p-1}\\
&\quad\times\sum_{a=1}^p\abs{
\skakko{{\rm E}\prod_{j\in\nu_a}X_{i_{(j)}}(s_{(j)})}-
\skakko{{\rm E}\prod_{j\in\nu_a \cap \upsilon}X_{i_{(j)}}(s_{(j)})}
\skakko{{\rm E}\prod_{j\in\nu_a \cap\upsilon^\mathsf{c}}X_{i_{(j)}}^*(s_{(j)}}}
\\
\leq&
8 \alpha^{\frac{1}{2}}(s_{(j^*+1)}-s_{(j^*)})\sum_{(\nu_1,\ldots,\nu_p)}(p-1)!
C_{\rm ineq,\ell}^{p-1}\\
&\quad\times
\sum_{a=1}^p
\skakko{{\rm E}\prod_{j\in\nu_a \cap \upsilon}X^{4}_{i_{(j)}}(s_{(j)})}^{\frac{1}{4}} 
\skakko{{\rm E}\prod_{j\in\nu_a \cap \upsilon^\mathsf{c}}X^{4}_{i_{(j)}}(s_{(j)})}^{\frac{1}{4}}\\
\leq&
C^\prime_{\ell}
\rho^{\frac{s_{(j^*+1)}-s_{(j^*)}}{2}}\\
\leq&
C^\prime_{\ell}
\rho^{\frac{\sum_{j=1}^{\ell-1}(s_{(j+1)}-s_{(j)})}{2(\ell-1)}}
\end{align*}
where $X_{i}^*(s)$ denotes an independent copy of $X_{i}(s)$, the summation $\sum_{(\nu_1,\ldots,\nu_p)}$ extends over all partitions $(\nu_1,\ldots,\nu_p)$ of $\{1,2,\cdots,\ell\}$,
\begin{align*}
&C_{\rm ineq, \ell}:=
\sup_{\substack{j\in\{1,\ldots,\ell\} \\ i_1,\ldots,i_j\in\{1,\ldots,K\} \\ (s_2,\ldots,s_j)\in\mathbb Z}}
\abs{{\rm E}X_{i_1}(0)X_{i_2}(s_2)\cdots X_{i_j}(s_j)},
\end{align*}
and
\begin{align*}
&C^\prime_{\ell}:=8C_\alpha\sup_{\substack{i_1,\ldots,i_\ell\in\{1,\ldots,K\} \\ (s_2,\ldots,s_\ell)\in\mathbb Z}}\Bigg|\sum_{(\nu_1,\ldots,\nu_p)}(p-1)!
C_{\rm ineq,\ell}^{p-1}
\\
&
\times\sum_{a=1}^p
\skakko{{\rm E}\prod_{j\in\nu_a \cap \upsilon}X^{4}_{i_{(j)}}(s_{(j)})}^{\frac{1}{4}} 
\skakko{{\rm E}\prod_{j\in\nu_a \cap \upsilon^\mathsf{c}}X^{4}_{i_{(j)}}(s_{(j)})}^{\frac{1}{4}}\Bigg|
\end{align*}
with the convention that ${\rm E}\skakko{\prod_{j\in\emptyset}X_{i_{(j)}}(s_{(j)})}=1$.
Therefore, for each $\ell\in\mathbb N$, there exist constants $C_\ell^{\prime\prime}\geq1$ and $\rho^\prime\in(0,1)$ such that 
\begin{align}\label{cum_bound}
\abs{{\rm cum}\{X_{i_1}(0),X_{i_2}(s_2),\ldots, X_{i_\ell}(s_\ell)\}}\leq C_\ell^{\prime\prime}\skakko{\rho^\prime}^{\sum_{j=1}^{\ell-1}(s_{(j+1)}-s_{(j)})}
\end{align}
for any $s_2,\ldots,s_{\ell}\in\mathbb Z$ and any $i_1,\ldots,i_\ell\in\{1,\ldots,K\}$.

Next, we show the summability condition for the cumulant (3) by using \eqref{cum_bound}. Simple algebra gives
\begin{align*}
&\sum_{s_2,\ldots,s_\ell=-\infty}^\infty\skakko{1+\sum_{j=2}^\ell\abs{s_j}^d}\abs{{\rm cum}\{X_{i_1}(0),X_{i_2}(s_2),\ldots, X_{i_\ell}(s_\ell)\}}\\
\leq&
C_\ell^{\prime\prime}
\sum_{s_2,\ldots,s_\ell=-\infty}^\infty
\skakko{1+\sum_{j=2}^\ell\abs{s_j}^d}
\skakko{\rho^\prime}^{\sum_{j=1}^{\ell-1}(s_{(j+1)}-s_{(j)})}\\
\leq&
C_\ell^{\prime\prime}(\ell-1)!\sum_{I}
\skakko{1+\sum_{j=2}^\ell\abs{s_j}^d}
\skakko{\rho^\prime}^{\sum_{j=1}^{\ell-1}(s_{(j+1)}-s_{(j)})}\\
\leq&
C_\ell^{\prime\prime}(\ell-1)!
\Bigg(
\sum_{s_{\ell}=-\infty}^{0}
\sum_{s_{\ell-1}=s_{\ell}}^{0}
\sum_{s_{\ell-2}=s_{\ell-1}}^{0}
\dots
\sum_{s_{3}=s_{4}}^{0}
\sum_{s_{2}=s_{3}}^{0}
\skakko{\sum_{j=2}^{\ell}|s_j|^d}
\skakko{\rho^\prime}^{-s_{\ell}}\\
&+
\sum_{\tau=3}^{\ell}
\sum_{s_{\ell}=-\infty}^{0}
\sum_{s_{\ell-1}=s_{\ell}}^{0}
\sum_{s_{\ell-2}=s_{\ell-1}}^{0}
\dots
\sum_{s_{\tau+1}=s_{\tau+2}}^{0}
\sum_{s_{\tau}=s_{\tau+1}}^{0}
\sum_{s_{2}=0}^{\infty}
\sum_{s_{3}=0}^{s_{2}}
\sum_{s_{4}=0}^{s_{3}}
\dots
\sum_{s_{\tau-2}=0}^{s_{\tau-3}}
\sum_{s_{\tau-1}=0}^{s_{\tau-2}}\\
&\skakko{\sum_{j=2}^{\ell}|s_j|^d}
\skakko{\rho^\prime}^{s_{2}-s_{\ell}}\\
&+
\sum_{s_{2}=0}^{\infty}
\sum_{s_{3}=0}^{s_{2}}
\sum_{s_{4}=0}^{s_{3}}
\dots
\sum_{s_{\ell-1}=0}^{s_{\ell-2}}
\sum_{s_{\ell}=0}^{s_{\ell-1}}
\skakko{\sum_{j=2}^{\ell}|s_j|^d}
\skakko{\rho^\prime}^{s_{2}}
\Bigg)\\
=&C_\ell^{\prime\prime}(\ell-1)!(L_1+L_2+L_3)
\end{align*}
where $\sum_{I}$ runs over all integers $s_2,\ldots,s_{\ell}$ such that $-\infty<s_{\ell}\leq s_{\ell-1}\leq\ldots\leq s_2<\infty$, 
\begin{align*}
L_1:=&\sum_{s_{\ell}=-\infty}^{0}
\sum_{s_{\ell-1}=s_{\ell}}^{0}
\sum_{s_{\ell-2}=s_{\ell-1}}^{0}
\dots
\sum_{s_{3}=s_{4}}^{0}
\sum_{s_{2}=s_{3}}^{0}
\skakko{\sum_{j=2}^{\ell}|s_j|^d}
\skakko{\rho^\prime}^{-s_{\ell}}
,\\
L_2:=&\sum_{\tau=3}^{\ell}
\sum_{s_{\ell}=-\infty}^{0}
\sum_{s_{\ell-1}=s_{\ell}}^{0}
\sum_{s_{\ell-2}=s_{\ell-1}}^{0}
\dots
\sum_{s_{\tau+1}=s_{\tau+2}}^{0}
\sum_{s_{\tau}=s_{\tau+1}}^{0}
\sum_{s_{2}=0}^{\infty}
\sum_{s_{3}=0}^{s_{2}}
\sum_{s_{4}=0}^{s_{3}}
\dots
\sum_{s_{\tau-2}=0}^{s_{\tau-3}}
\sum_{s_{\tau-1}=0}^{s_{\tau-2}}\\
&\skakko{\sum_{j=2}^{\ell}|s_j|^d}
\skakko{\rho^\prime}^{s_{2}-s_{\ell}}
,\\
\text{and }
L_3:=&\sum_{s_{2}=0}^{\infty}
\sum_{s_{3}=0}^{s_{2}}
\sum_{s_{4}=0}^{s_{3}}
\dots
\sum_{s_{\ell-1}=0}^{s_{\ell-2}}
\sum_{s_{\ell}=0}^{s_{\ell-1}}
\skakko{\sum_{j=2}^{\ell}|s_j|^d}
\skakko{\rho^\prime}^{s_{2}}
.
\end{align*}
We only prove $L_3<\infty$, as the proofs of $L_1,L_2<\infty$ is similar. It can be seen that
\begin{align*}
L_3\leq&
(\ell-1)
\sum_{s_{2}=0}^{\infty}
|s_2|^d\skakko{\rho^\prime}^{s_{2}}
\sum_{s_{3}=0}^{s_{2}}
\sum_{s_{4}=0}^{s_{3}}
\dots
\sum_{s_{\ell-1}=0}^{s_{\ell-2}}
\sum_{s_{\ell}=0}^{s_{\ell-1}}
1\\
\leq&
(\ell-1)
\sum_{s_{2}=0}^{\infty}
|s_2|^d\skakko{\rho^\prime}^{s_{2}}
\sum_{s_{3}=0}^{s_{2}}
\sum_{s_{4}=0}^{s_{2}}
\dots
\sum_{s_{\ell-1}=0}^{s_{2}}
\sum_{s_{\ell}=0}^{s_{2}}
1\\
\leq&
(\ell-1)
\sum_{s_{2}=0}^{\infty}
|s_2|^{d+\ell-2}\skakko{\rho^\prime}^{s_{2}}
\end{align*}
which is finite. \qed

\subsection{Proof of Theorem 2.1}

First, we show the uniqueness of the processes $G_1$ and $G_2$. 
Suppose that there exists the processes $\{G_d(t)\}$ for $d=1,\ldots,K$, which taking the form of
\begin{align*}
G_d(t):=\sum_{j=1}^d\sum_{k_j=-\infty}^\infty a_{dj}(k_j)X_j(t-k_j)
\end{align*}
such that $X_0(t)=\zeta + \sum_{d=1}^KG_d(t)+\epsilon(t)$, $G_i$ and $G_j$ are orthogonal for any $i,j(\neq i)\in\{1,\ldots,K\}$, i.e., ${\rm E}G_i(t)G_j(t^\prime)=0$ for any $t,t^\prime\in\mathbb Z$, $\sum_{k=-\infty}^\infty|a_{dj}(k)|<\infty$ for any $d\in \{1,\ldots,K\}$ and $j\in\{1,\ldots,d\}$, and the transfer function $A_{dd}\skakko{e^{-\mathrm{i}\lambda}}:=\sum_{k=-\infty}^\infty a_{dd}(k)e^{-\mathrm{i}k\lambda}$ satisfies $A_{dd}\skakko{e^{-\mathrm{i}\lambda}}\neq0$ for all $\lambda\in[-\pi,\pi]$ and $d\in\{1,\ldots,K\}$.
From the orthogonality ${\rm E}G_d(t)G_{d^\prime}(t^\prime)=0$ for $d\in\{2,\ldots,K\}$, $d^\prime\in\{1,\ldots,d-1\}$, and all $t,t^\prime\in\mathbb Z$,
it holds that
\begin{align}\nonumber
&f_{G_dG_{d^\prime}}(\lambda)=0\\\nonumber
\Leftrightarrow
&\sum_{j=1}^d\sum_{j^\prime=1}^{d^\prime}
A_{dj} \skakko{e^{-i\lambda}}A_{d^\prime j^\prime}  \skakko{e^{i\lambda}}f_{jj^\prime}(\lambda)=0\\\nonumber
\Leftrightarrow
&A_{d^\prime d^\prime}  \skakko{e^{i\lambda}}\sum_{j=1}^d
A_{dj} \skakko{e^{-i\lambda}}f_{jd^\prime}(\lambda)=0\\\label{fGG_eq}
\Leftrightarrow
&\sum_{j=1}^d
A_{dj} \skakko{e^{-i\lambda}}f_{jd^\prime}(\lambda)=0,\quad 
\end{align}
and thus
\begin{align*}
\bm{f}_{d-1}^\top(\lambda){\bm{A}_{d}^\flat} \skakko{e^{-i\lambda}}
=-A_{dd} \skakko{e^{-i\lambda}}\overline{\bm{f}_{d-1,d}^\flat(\lambda)},
\end{align*}
where
\begin{align*}
{\bm{A}_{d}^\flat} \skakko{e^{-i\lambda}}:=\skakko{{A}_{d1} \skakko{e^{-i\lambda}},\ldots,{A}_{d(d-1)} \skakko{e^{-i\lambda}}}^\top\quad\text{and }
\bm{f}_{a,b}^\flat:=&({f}_{1b}(\lambda),\ldots,{f}_{ab}(\lambda))^\top
.
\end{align*}
Since $\bm{f}_K(\lambda)$ is positive definite if and only if all principal sub-matrices is positive definite, we obtain
\begin{align}\label{Add_sol1}
{\bm{A}_{d}^\flat} \skakko{e^{-i\lambda}}
=-A_{dd} \skakko{e^{-i\lambda}}{\bm{f}_{d-1}^\top}^{-1}(\lambda)\overline{\bm{f}_{d-1,d}^\flat(\lambda)},
\end{align}
or equivalently by Cramer's rule, 
\begin{align}\label{Add_sol1_2}
A_{dj} \skakko{e^{-i\lambda}}
=-A_{dd} \skakko{e^{-i\lambda}}\frac{{\rm det}\skakko{\overline{\bm{f}_{j,d}^\ddag(\lambda)}}}{{\rm det}\skakko{\bm{f}_{d-1}^\top(\lambda)}}\quad\text{for $j=1,\ldots,d-1$},
\end{align}
where
\begin{align*}
{\bm{f}}_{j,d}^\ddag(\lambda):=\skakko{{\bm{f}_{d-1,1}^\flat(\lambda)},\ldots,{\bm{f}_{d-1,j-1}^\flat(\lambda)},{\bm{f}_{d-1,d}^\flat(\lambda)},{\bm{f}_{d-1,j+1}^\flat(\lambda)},\ldots,{\bm{f}_{d-1,d-1}^\flat(\lambda)}}.
\end{align*}
From  \eqref{Add_sol1}, \eqref{Add_sol1_2}, and 
\begin{align*}
\bm{f}_{d}(\lambda)=
\begin{pmatrix}
\bm{f}_{d-1}(\lambda)&\bm{f}_{d-1,d}^\flat(\lambda)\\
{\bm{f}_{d-1,d}^\flat}^*(\lambda)&{f}_{dd}(\lambda) \\
\end{pmatrix},
\end{align*}
simple algebra yields
\begin{align}\nonumber
{\rm det}\skakko{\bm{f}_{d}(\lambda)}=&
{\rm det}\skakko{\bm{f}_{d-1}(\lambda)}
{\rm det}\skakko{{f}_{dd}(\lambda) - {\bm{f}_{d-1,d}^\flat}^*(\lambda){\bm{f}_{d-1}^{-1}}(\lambda)\bm{f}_{d-1,d}^\flat(\lambda)}\\\nonumber
=&{\rm det}\skakko{\bm{f}_{d-1}(\lambda)}
{\rm det}\skakko{{f}_{dd}(\lambda) - {\bm{f}_{d-1,d}^\flat}^\top(\lambda){{\bm{f}_{d-1}^\top}^{-1}}(\lambda)\overline{\bm{f}_{d-1,d}^\flat(\lambda)}}\\\nonumber
=&{\rm det}\skakko{\bm{f}_{d-1}(\lambda)}
{\rm det}\skakko{{f}_{dd}(\lambda) + \skakko{A_{dd} \skakko{e^{-i\lambda}}}^{-1}{\bm{f}_{d-1,d}^\flat}^\top(\lambda){\bm{A}_{d}^\flat} \skakko{e^{-i\lambda}}}\\\nonumber
=&{\rm det}\skakko{\bm{f}_{d-1}(\lambda)}
{\rm det}\skakko{{f}_{dd}(\lambda) + \skakko{A_{dd} \skakko{e^{-i\lambda}}}^{-1}
\sum_{j=1}^{d-1}{f}_{jd}(\lambda){A}_{dj} \skakko{e^{-i\lambda}}
}\\\nonumber
=&{\rm det}\skakko{\bm{f}_{d-1}(\lambda)}
{\rm det}\skakko{{f}_{dd}(\lambda) - 
\sum_{j=1}^{d-1}{f}_{jd}(\lambda)\frac{{\rm det}\skakko{\overline{\bm{f}_{j,d}^\ddag(\lambda)}}}{{\rm det}\skakko{\bm{f}_{d-1}^\top(\lambda)}}
}\\\label{f_formula}
=&
{f}_{dd}(\lambda){\rm det}\skakko{\bm{f}_{d-1}(\lambda)}
-
\sum_{j=1}^{d-1}{f}_{jd}(\lambda){\rm det}\skakko{\overline{\bm{f}_{j,d}^\ddag(\lambda)}}.
\end{align}
By the orthogonality ${\rm E}G_d(t)\overline{X_{0}(t^\prime)}={\rm E}G_d(t)\overline{G_d(t^\prime)}$ for $d\in\{2,\ldots,K\}$ and all $t,t^\prime\in\mathbb Z$, \eqref{fGG_eq}, \eqref{Add_sol1_2}, \eqref{f_formula}, 
we obtain
\begin{align*}
&f_{G_{d}0}(\lambda)=f_{G_dG_d}(\lambda)\\
\Leftrightarrow&
\sum_{j=1}^dA_{dj} \skakko{e^{-i\lambda}}f_{j0}(\lambda)=
\sum_{j,j^\prime=1}^dA_{dj} \skakko{e^{-i\lambda}}A_{dj^\prime}  \skakko{e^{i\lambda}}f_{jj^\prime}(\lambda)\\
\Leftrightarrow&
\sum_{j=1}^dA_{dj} \skakko{e^{-i\lambda}}f_{j0}(\lambda)=
\sum_{j^\prime=1}^dA_{dj^\prime}  \skakko{e^{i\lambda}}
\sum_{j=1}^dA_{dj} \skakko{e^{-i\lambda}}f_{jj^\prime}(\lambda)
\\
\Leftrightarrow&
\sum_{j=1}^dA_{dj} \skakko{e^{-i\lambda}}f_{j0}(\lambda)=
A_{dd}  \skakko{e^{i\lambda}}
\sum_{j=1}^dA_{dj} \skakko{e^{-i\lambda}}f_{jd}(\lambda)
\\
\Leftrightarrow&
\sum_{j=1}^{d-1}A_{dj} \skakko{e^{-i\lambda}}f_{j0}(\lambda)+A_{dd} \skakko{e^{-i\lambda}}f_{d0}(\lambda)\\\nonumber
&=
A_{dd}  \skakko{e^{i\lambda}}\skakko{
\sum_{j=1}^{d-1}A_{dj} \skakko{e^{-i\lambda}}f_{jd}(\lambda)+A_{dd} \skakko{e^{-i\lambda}}f_{dd}(\lambda)}
\\
\Leftrightarrow&
-\sum_{j=1}^{d-1}{\rm det}\skakko{\overline{\bm{f}_{j,d}^\ddag(\lambda)}}f_{j0}(\lambda)+f_{d0}(\lambda){\rm det}\skakko{\bm{f}_{d-1}^\top(\lambda)}\\
&=
A_{dd}  \skakko{e^{i\lambda}}\skakko{
-\sum_{j=1}^{d-1}{\rm det}\skakko{\overline{\bm{f}_{j,d}^\ddag(\lambda)}}f_{jd}(\lambda)+f_{dd}(\lambda){\rm det}\skakko{\bm{f}_{d-1}^\top(\lambda)}}
\\
\Leftrightarrow&
A_{dd}  \skakko{e^{i\lambda}}
=
\frac{-\sum_{j=1}^{d-1}{\rm det}\skakko{\overline{\bm{f}_{j,d}^\ddag(\lambda)}}f_{j0}(\lambda)+f_{d0}(\lambda){\rm det}\skakko{\bm{f}_{d-1}^\top(\lambda)}}{-\sum_{j=1}^{d-1}{\rm det}\skakko{\overline{\bm{f}_{j,d}^\ddag(\lambda)}}f_{jd}(\lambda)+f_{dd}(\lambda){\rm det}\skakko{\bm{f}_{d-1}^\top(\lambda)}}
\\
\Leftrightarrow&
A_{dd}  \skakko{e^{i\lambda}}
=
\frac{-\sum_{j=1}^{d-1}{\rm det}\skakko{\overline{\bm{f}_{j,d}^\ddag(\lambda)}}f_{j0}(\lambda)+f_{d0}(\lambda){\rm det}\skakko{\bm{f}_{d-1}^\top(\lambda)}}{{\rm det}\skakko{\bm{f}_{d}(\lambda)}}
\\
\end{align*}
which also gives, for $d=2,\ldots,K$ and $d^\prime\in\{1,\ldots,d-1\}$,
\begin{align*}
A_{dd^\prime} \skakko{e^{i\lambda}}
=&-A_{dd} \skakko{e^{i\lambda}}\frac{{\rm det}\skakko{\bm{f}_{d^\prime,d}^\ddag(\lambda)}}{{\rm det}\skakko{\bm{f}_{d-1}(\lambda)}}\\
\text{ and }
f_{G_dG_d}(\lambda)
=&
A_{dd}  \skakko{e^{i\lambda}}
\sum_{j=1}^dA_{dj} \skakko{e^{-i\lambda}}f_{jd}(\lambda)\\
=&
A_{dd}  \skakko{e^{i\lambda}}\skakko{
\sum_{j=1}^{d-1}A_{dj} \skakko{e^{-i\lambda}}f_{jd}(\lambda)
+A_{dd} \skakko{e^{-i\lambda}}f_{dd}(\lambda)}\\
=&
\frac{\left|A_{dd}  \skakko{e^{i\lambda}}\right|^2}{{\rm det}\skakko{\bm{f}_{d-1}^\top(\lambda)}}\skakko{
-\sum_{j=1}^{d-1}{\rm det}\skakko{\overline{\bm{f}_{j,d}^\ddag(\lambda)}}f_{jd}(\lambda)+f_{dd}(\lambda){\rm det}\skakko{\bm{f}_{d-1}^\top(\lambda)}}\\
=&
\frac{\left|A_{dd}  \skakko{e^{i\lambda}}\right|^2}{{\rm det}\skakko{\bm{f}_{d-1}(\lambda)}}{\rm det}\skakko{\bm{f}_{d}(\lambda)}.
\end{align*}
For the case $d=1$, we use the orthogonality ${\rm E}G_1(t)\overline{X_{0}(t^\prime)}={\rm E}G_1(t)\overline{G_1(t^\prime)}$ for all $t,t^\prime\in\mathbb Z$ and have 
\begin{align*}
f_{G_{1}0}(\lambda)=f_{G_1G_1}(\lambda)
\Leftrightarrow&
A_{11} \skakko{e^{-i\lambda}}f_{10}(\lambda)=
\left|A_{11} \skakko{e^{-i\lambda}}\right|^2f_{11}(\lambda)\\
\Leftrightarrow&
A_{11} \skakko{e^{i\lambda}}
=\frac{f_{10}(\lambda)}{f_{11}(\lambda)}
\end{align*}
and 
\begin{align*}
f_{G_1G_1}(\lambda)=
\frac{\left|f_{10}(\lambda)\right|^2}{f_{11}(\lambda)}.
\end{align*}

To complete the proof, we need show that,
$A_{KK}\skakko{e^{-\mathrm{i}\lambda }}=0$ for all $\lambda\in[-\pi,\pi]$ if and only if 
\begin{align}\label{AKKright}
\frac{-\sum_{i=1}^{K-1}{\rm det}\skakko{\overline{\bm{f}_{i,K}^\ddag(\lambda)}}f_{i0}(\lambda)+{\rm det}\skakko{\bm{f}_{K-1}(\lambda)}f_{K0}(\lambda)}{{\rm det}\skakko{\bm{f}_{K}(\lambda)}}=0 \quad \text{for all $\lambda\in[-\pi,\pi]$}.
\end{align} 
First, we prove that
$A_{KK}\skakko{e^{-\mathrm{i}\lambda }}=0$ for all $\lambda\in[-\pi,\pi]$ is necessary for 
$$
\frac{-\sum_{i=1}^{K-1}{\rm det}\skakko{\overline{\bm{f}_{i,K}^\ddag(\lambda)}}f_{i0}(\lambda)+{\rm det}\skakko{\bm{f}_{K-1}(\lambda)}f_{K0}(\lambda)}{{\rm det}\skakko{\bm{f}_{K}(\lambda)}}=0 \quad \text{for all $\lambda\in[-\pi,\pi]$}
$$
by contraposition. Suppose $A_{KK}\skakko{e^{-\mathrm{i}\lambda }}\neq0$ for some $\lambda\in[-\pi,\pi]$, then from the above discussion 
 we have
$$A_{KK}\skakko{e^{-\mathrm{i}\lambda }}=
\frac{-\sum_{i=1}^{K-1}{\rm det}\skakko{\overline{\bm{f}_{i,K}^\ddag(\lambda)}}f_{i0}(\lambda)+{\rm det}\skakko{\bm{f}_{K-1}(\lambda)}f_{K0}(\lambda)}{{\rm det}\skakko{\bm{f}_{K}(\lambda)}}
$$
which is not zero for some $\lambda\in[-\pi,\pi]$. 
Next, we show
$A_{KK}\skakko{e^{-\mathrm{i}\lambda }}=0$ for all $\lambda\in[-\pi,\pi]$ is sufficient for 
$$
\frac{-\sum_{i=1}^{K-1}{\rm det}\skakko{\overline{\bm{f}_{i,K}^\ddag(\lambda)}}f_{i0}(\lambda)+{\rm det}\skakko{\bm{f}_{K-1}(\lambda)}f_{K0}(\lambda)}{{\rm det}\skakko{\bm{f}_{K}(\lambda)}}=0 \quad \text{for all $\lambda\in[-\pi,\pi]$}.
$$
The condition $A_{KK}\skakko{e^{-\mathrm{i}\lambda }}=0$ for all $\lambda\in[-\pi,\pi]$ implies that 
$A_{K,j}(e^{-i\lambda})=0$ for $j=1,\ldots,K-1$ and $G_{d}(t)=0$ for all $t\in\mathbb Z$. Further, for $j=1,\dots,d$, and any $t,t^\prime\in\mathbb Z$, $EX_{0}(t)\overline{X_{j}(t^{\prime})}=\sum_{i=1}^{d-1}EG_{i}(t)\overline{X_{j}(t^{\prime})}$ yields 
\begin{equation*}
	f_{0,j}(\lambda)=\sum_{i=1}^{d-1}\sum_{k=i}^{d-1}A_{k,i}(e^{-i\lambda})f_{i,j}(\lambda),
\end{equation*}
which gives that
\begin{align*}
		&-\sum_{j=1}^{d-1}\det(\bm{f}_{j,d}^{\ddag}(\lambda))f_{0,j}(\lambda)+\det(\bm{f}_{d-1}(\lambda))f_{0,d}(\lambda)\\
	=&-\sum_{j=1}^{d-1}\det(\bm{f}_{j,d}^{\ddag}(\lambda))\sum_{i=1}^{d-1}\sum_{k=i}^{d-1}A_{k,i}(e^{-i\lambda})f_{i,j}(\lambda)+\det(\bm{f}_{d-1}(\lambda))\sum_{i=1}^{d-1}\sum_{k=i}^{d-1}A_{k,i}(e^{-i\lambda})f_{i,d}(\lambda)\\
	=&\sum_{i=1}^{d-1}\sum_{k=i}^{d-1}A_{k,i}(e^{-i\lambda})\left(-\sum_{j=1}^{d-1}\det(\bm{f}_{j,d}^{\ddag}(\lambda))f_{i,j}(\lambda)+\det(\bm{f}_{d-1}(\lambda))f_{i,d}(\lambda)\right).
\end{align*}
Therefore, it is sufficient to show 
\begin{equation*}
-\sum_{j=1}^{d-1}\det(\bm{f}_{j,d}^{\ddag}(\lambda))f_{i,j}(\lambda)+\det(\bm{f}_{d-1}(\lambda))f_{i,d}(\lambda)=0
\end{equation*}
for $i=1,\dots,d-1$. Denote the minor of matrix $\bm{A}$ with $i$th row and $j$th column eliminated as $\bm{f}^{(i,j)}$ and note that $\det(\bm{f}_{j,d}^{\ddag(i,j)}(\lambda))= \det(\bm{f}_{d-1}^{(i,j)}(\lambda))$.
Then, 
\begin{align}\nonumber
&-\sum_{j=1}^{d-1}\det(\bm{f}_{j,d}^{\ddag}(\lambda))f_{i,j}(\lambda)+\det(\bm{f}_{d-1}(\lambda))f_{i,d}(\lambda)\\\nonumber
   =&-\sum_{j=1}^{d-1}
   \left(\sum_{h=1}^{d-1}(-1)^{j+h}\det(\bm{f}_{j,d}^{\ddag(h,j)}(\lambda))f_{h,d}(\lambda)\right)f_{i,j}(\lambda)\\\nonumber
&+\left(\sum_{j=1}^{d-1}(-1)^{j+i}\det(\bm{f}_{d-1}^{(i,j)}(\lambda))f_{i,j}(\lambda)\right)f_{i,d}(\lambda)\\\nonumber
   =&-\sum_{h=1}^{d-1}
   \left(\sum_{j=1}^{d-1}(-1)^{j+h}\det(\bm{f}_{j,d}^{\ddag(h,j)}(\lambda))f_{i,j}(\lambda)\right)f_{h,d}(\lambda)\\\nonumber
&+\left(\sum_{j=1}^{d-1}(-1)^{j+i}\det(\bm{f}_{d-1}^{(i,j)}(\lambda))f_{i,j}(\lambda)\right)f_{i,d}(\lambda)\\\nonumber
=&-\sum_{h=1,h\neq i}^{d-1}
		\left(\sum_{j=1}^{d-1}(-1)^{j+h}\det(\bm{f}_{j,d}^{\ddag(h,j)}(\lambda))f_{i,j}(\lambda)\right)f_{h,d}(\lambda)\\\label{eq:proof_thm2.1}
		=&-\sum_{h=1,h\neq i}^{d-1}
		\left(\sum_{j=1}^{d-1}(-1)^{j+h}\det(\bm{f}_{d-1}^{(h,j)}(\lambda))f_{i,j}(\lambda)\right)f_{h,d}(\lambda).
\end{align}
The determinant of the matrix
\begin{align*}
\begin{pmatrix}
f_{11}&\cdots&f_{1(d-1)}\\
\vdots&&\\
f_{(h-1)1}&\cdots&f_{(h-1)(d-1)}\\
f_{i1}&\cdots&f_{i(d-1)}\\
f_{(h+1)1}&\cdots&f_{(h+1)(d-1)}\\
\vdots&&\\
f_{(d-1)1}&\cdots&f_{(d-1)(d-1)}\\
\end{pmatrix}
\end{align*}
can be expanded as $\sum_{j=1}^{d-1}(-1)^{j+h}\det(\bm{f}_{d-1}^{(h,j)}(\lambda))f_{i,j}(\lambda)$ but is zero since the $i$-th row and $h$-th row of the matrix are the same, which yields that  \eqref{eq:proof_thm2.1} is zero.

Second, we show the existence of the processes $\{G_d(t)\}$. To this end, we shall confirm that $\{G_d(t)\}$ defined by (5) with $a_{dj}$ whose transfer function is given by (6) for $d=j=1$, (7) for $d=j\in\{2,\ldots,K\}$, (8) for $d\in\{2,\ldots,K\}$ and $j\in\{2,\ldots,d-1\}$ satisfies 
(i) $G_i$ and $G_j$ are orthogonal for any $i,j(\neq i)\in\{1,\ldots,K\}$, 
(ii) $X_0(t)=\zeta + \sum_{j=1}^KG_j(t)+\epsilon(t)$, 
(iii) $\sum_{k=-\infty}^\infty|a_{dj}(k)|<\infty$ for any $d\in \{1,\ldots,K\}$ and $j\in\{1,\ldots,d\}$, and 
(iv) $A_{dd}\skakko{e^{-\mathrm{i}\lambda}}\neq0$ for all $\lambda\in[-\pi,\pi]$ and $d\in\{1,\ldots,K\}$.

For $d\leq K$ and $j\leq d-1$, Cramer's rule tells us that 
$A_{dj}$ satisfies 
\begin{align*}
\bm{f}_{d-1}^\top(\lambda){\bm{A}_{d}^\flat} \skakko{e^{-i\lambda}}
=-A_{dd} \skakko{e^{-i\lambda}}\overline{\bm{f}_{d-1,d}^\flat(\lambda)},
\end{align*}
and thus we have, for $d^\prime=1,\ldots,d-1$,
\begin{align*}
&\sum_{j=1}^d
A_{dj} \skakko{e^{-i\lambda}}f_{jd^\prime}(\lambda)=0\\
\Leftrightarrow
&A_{d^\prime d^\prime}  \skakko{e^{i\lambda}}\sum_{j=1}^d
A_{dj} \skakko{e^{-i\lambda}}f_{jd^\prime}(\lambda)=0\\
\Leftrightarrow
&\sum_{j=1}^d\sum_{j^\prime=1}^{d^\prime}
A_{dj} \skakko{e^{-i\lambda}}A_{d^\prime j^\prime}  \skakko{e^{i\lambda}}f_{jj^\prime}(\lambda)=0\\
\Leftrightarrow
&f_{G_dG_{d^\prime}}(\lambda)=0,
\end{align*}
which shows (i).

In order to verify (ii), we show $B_j\skakko{e^{-i\lambda}}=\sum_{k=j}^dA_{k,j}\skakko{e^{-i\lambda}}$.
The orthogonality of $G_d$'s gives
$\sum_{j=1}^d
A_{dj} \skakko{e^{-i\lambda}}f_{jd^\prime}(\lambda)=0,
$ which is equivalent to ${\rm E}G_d(t)\bar Y_{d^\prime}(t^\prime)=0$ for any  $t,t^\prime\in\mathbb Z$, $d\leq K$, and $d^\prime=1,\ldots,d-1$. From ${\rm E}Y_0(t)\bar G_d(t^{\prime})={\rm E}G_d(t)\bar G_d(t^{\prime})$, it holds that
\begin{align}\label{eq_ii}
\sum_{j=1}^{d}A_{dj}\skakko{e^{i\lambda}}f_{0j}(\lambda)=&B_d\skakko{e^{-i\lambda}}\sum_{j=1}^{d}A_{dj}(e^{i\lambda})f_{dj}(\lambda).
\end{align}
By noting that 
${\rm det}\skakko{\bm{f}_{d}(\lambda)}=
-\sum_{i=1}^{d-1}{\rm det}\skakko{\overline{\bm{f}_{i,d}^\ddag(\lambda)}}f_{id}(\lambda)+{\rm det}\skakko{\bm{f}_{d-1}(\lambda)}f_{dd}(\lambda)
$, \eqref{eq_ii} yields that 
\begin{align*}
B_d\skakko{e^{-i\lambda}}:=
&
\frac{\sum_{j=1}^{d}A_{dj}\skakko{e^{i\lambda}}f_{0j}(\lambda)}{\sum_{j=1}^{d}A_{dj}(e^{i\lambda})f_{dj}(\lambda)}\\
=&
\frac{-\sum_{i=1}^{d-1}{\rm det}\skakko{\overline{\bm{f}_{i,d}^\ddag(\lambda)}}f_{i0}(\lambda)+{\rm det}\skakko{\bm{f}_{d-1}(\lambda)}f_{d0}(\lambda)}{-\sum_{i=1}^{d-1}{\rm det}\skakko{\overline{\bm{f}_{i,d}^\ddag(\lambda)}}f_{id}(\lambda)+{\rm det}\skakko{\bm{f}_{d-1}(\lambda)}f_{dd}(\lambda)}\\
=&
A_{dd}\skakko{e^{-i\lambda}}.
\end{align*}
Consider ${\rm E}Y_0(t)\bar G_{d-1}(t^{\prime})={\rm E}G_{d-1}(t)\bar G_{d-1}(t^{\prime})$, then we obtain
\begin{align*}
		&\sum_{j=1}^{d-1}A_{(d-1)j}\skakko{e^{i\lambda}}f_{0j}(\lambda)\\
		=&B_{d-1}\skakko{e^{-i\lambda}}\sum_{j=1}^{d-1}A_{(d-1)j}\skakko{e^{i\lambda}}f_{d-1j}(\lambda)+A_{dd}\skakko{e^{-i\lambda}}\sum_{j=1}^{d-1}A_{(d-1)j}\skakko{e^{i\lambda}}f_{dj}(\lambda).
\end{align*}
Somewhat lengthy calculations gives
\begin{align*}
&B_{d-1}\skakko{e^{-i\lambda}}\\
=&\frac{\sum_{j=1}^{d-1}A_{(d-1)j}\skakko{e^{i\lambda}}f_{0j}(\lambda)}{\sum_{j=1}^{d-1}A_{(d-1)j}\skakko{e^{i\lambda}}f_{(d-1)j}(\lambda)}-A_{dd}\skakko{e^{-i\lambda}}\frac{\sum_{j=1}^{d-1}A_{(d-1)j}\skakko{e^{i\lambda}}f_{dj}(\lambda)}{\sum_{j=1}^{d-1}A_{(d-1)j}\skakko{e^{i\lambda}}f_{(d-1)j}(\lambda)}\\
=&A_{d-1,d-1}\skakko{e^{-i\lambda}}+A_{d,d-1}\skakko{e^{-i\lambda}}.
\end{align*}
In the same manner, 
${\rm E}(Y_0(t)-G_d(t))\bar G_{d-1}(t^{\prime})$ and ${\rm E}(Y_0(t)-G_d(t))\bar G_{d-2}(t^{\prime})$ imply
\begin{align*}
B_{d-1}^*\skakko{e^{-i\lambda}}=&A_{(d-1)(d-1)}\skakko{e^{-i\lambda}}\\
B_{d-2}^*\skakko{e^{-i\lambda}}=&A_{(d-2)(d-2)}\skakko{e^{-i\lambda}}+A_{(d-1)(d-2)}\skakko{e^{-i\lambda}},
\end{align*}
respectively, where $B_{j}^{*}\skakko{e^{-i\lambda}}:=B_j\skakko{e^{-i\lambda}}-A_{dj}\skakko{e^{-i\lambda}}$. Moreover, 
${\rm E}(Y_0(t)-G_d(t)-G_{d-1}(t))\bar G_{d-2}(t^{\prime})$ and ${\rm E}(Y_0(t)-G_d(t)-G_{d-1}(t))\bar G_{d-3}(t^{\prime})$ imply
\begin{align*}
B_{d-1}^{**}\skakko{e^{-i\lambda}}=&A_{(d-2)(d-2)}\skakko{e^{-i\lambda}}\\
B_{d-2}^{**}\skakko{e^{-i\lambda}}=&A_{(d-3)(d-3)}\skakko{e^{-i\lambda}}+A_{(d-2)(d-3)}\skakko{e^{-i\lambda}},
\end{align*}
respectively, where $B_{j}^{**}\skakko{e^{-i\lambda}}:=B_j^{*}\skakko{e^{-i\lambda}}-A_{(d-1)j}\skakko{e^{-i\lambda}}$. 
By repeating the argument, we obtain $B_j\skakko{e^{-i\lambda}}=\sum_{k=j}^dA_{k,j}\skakko{e^{-i\lambda}}$.

The summability of linear filters (iii) follows from the twice continuously differentiability of $A_{ij}$ (see \citealt[Section 4.4, p.25]{katznelson04}).

The condition (iv) follows directly from the assumption of Theorem 2.1.
\qed

\subsection{Proof of Theorem 3.1}
Since $T_n$ and $\hat \mu_{n,K}$ are the functional of $\hat {\bm{f}}(\lambda)$,  write $T_n$ and $\hat \mu_{n,K}$ as $T\skakko{\hat {\bm{f}}(\lambda)}$ and $\mu\skakko{\hat {\bm{f}}(\lambda)}$, respectively, First, we show that the effect on the estimation of ${\rm E}X_i^\prime$ in the kernel density estimator is asymptotically negligible, that is,
\begin{align}\label{diff_Tf}
T\skakko{\hat {\bm{f}}(\lambda)}-T\skakko{\hat{\hat{\bm{f}}}(\lambda)}=o_p(1)\quad n\to\infty,
\end{align}
where 
\begin{align*}
\hat{\hat{f}}_{ij}(\lambda):=\frac{1}{2\pi}\sum_{h =1-n}^{n-1}\omega\skakko{\frac{h}{M_n}}\hat{\hat{{\gamma}}}_{ij}(h)e^{-\mathrm{i}h\lambda}
\end{align*}
with, for $h\in\{0,\ldots,n-1\}$,
$$\hat{\hat{\gamma}}_{ij}(h):=\frac{1}{n-h}\sum_{t=1}^{n-h}(X_i{(t+h)}- {\rm E}X_i{(t+h)})({ X}_j({t})- {\rm E}X_j{(t)}),$$
 for $h\in\{-n+1,\ldots,-1\}$,
$$\hat{\hat{\gamma}}_{ij}(h):=\frac{1}{n+{h}}\sum_{t=-h+1}^{n}(X_i{(t+h)}- {\rm E}X_i{(t+h)})({ X}_j({t})- {\rm E}X_j{(t)}).$$
From Assumption 3.1 (A2) and the fact that $\left|\hat{\hat{\gamma}}_{ij}(h)- {\hat{\gamma}}_{ij}(h)\right|=O_p(1/ n)$ uniformly in $h\in\{-n+1,\ldots,n-1\}$ (see, e.g., \citealt[Remark 3.1]{politis11}),  we obtain
\begin{align}\nonumber
\sup_{\lambda\in[-\pi,\pi]}\left|{\hat{f}}_{ij}(\lambda)- \hat{\hat{f}}_{ij}(\lambda)\right|
\leq&\frac{1}{2\pi}\sum_{h =1-n}^{n-1}\omega\skakko{\frac{h}{M_n}}\left|\hat{\hat{\gamma}}_{ij}(h)- {\hat{\gamma}}_{ij}(h)\right|\\\label{diff_f}
=&O_p\skakko{\frac{M_n}{n}}.
\end{align}
Then, Lipschitz continuity of $\mu$ gives
\begin{align*}
&{\rm P}\skakko{M_n\left|\mu\skakko{\hat {\bm{f}}(\lambda)}-\mu\skakko{\hat{\hat {\bm{f}}}{(\lambda)}}\right|>\epsilon}\\
\leq&
{\rm P}\skakko{M_n L_\mu
\int_{-\pi}^\pi
\left\| \hat {\bm{f}}(\lambda)-\hat{\hat {\bm{f}}}(\lambda)\right\|_2
{\rm d}\lambda
>\epsilon, \sup_{\lambda\in[-\pi,\pi]}\left\|\hat {\bm{f}}(\lambda)-{\bm{f}}(\lambda)\right\|_2 \leq\delta}\\
&+
{\rm P}\skakko{
\sup_{\lambda\in[-\pi,\pi]}\left\|\hat {\bm{f}}(\lambda)-{\bm{f}}(\lambda)\right\|_2>\delta}\\
\leq&
{\rm P}\skakko{M_n L_\mu
\int_{-\pi}^\pi
\left\| \hat {\bm{f}}(\lambda)-\hat{\hat {\bm{f}}}(\lambda)\right\|_2
{\rm d}\lambda
>\epsilon}
+
{\rm P}\skakko{
\sup_{\lambda\in[-\pi,\pi]}\left\|\hat {\bm{f}}(\lambda)-{\bm{f}}(\lambda)\right\|_2>\delta},
\end{align*}
where $L_\mu$ is the Lipschitz constant. 
\citet[Theorem 2.1]{robinson91} and \eqref{diff_f} yield that 
\begin{align}\label{diif_mu}
M_n\skakko{\mu\skakko{\hat {\bm{f}}(\lambda)}-\mu\skakko{\hat{\hat {\bm{f}}}{(\lambda)}}}=o_p(1)\quad \text{as $n\to\infty$.}
\end{align}
In addition, Theorems 5.9.1 and 7.4.1--7.4.4 of \cite{brillinger1981} yield that
\begin{align*}
{\hat{f}}_{ij}(\lambda)-{\check{f}}_{ij}(\lambda) =& O_p\skakko{\frac{M_n^2}{n}+\frac{M_n}{n}\log n}\quad \text{uniformly in $\lambda\in[-\pi,\pi]$}\\
\text{and }{\check{f}}_{ij}(\lambda)-{f}_{ij}(\lambda)=&O_p\skakko{\sqrt\frac{M_n}{n}+\frac{1}{M_n^2}}\quad \text{uniformly in $\lambda\in[-\pi,\pi]$},
\end{align*}
where ${\check{f}}_{ij}(\lambda)$ is  the discretized version of the smoothed periodogram defined as 
\begin{align*}
{\check{f}}_{ij}(\lambda):=\frac{2\pi}{n-u}\sum_{s=1}^{n-u-1}W^{(n)}\skakko{\lambda-\frac{2\pi s}{n-u}}I_{ij}\skakko{\frac{2\pi s}{n-u}}
\end{align*}
with $W^{(n)}(\lambda):=M_n\sum_{j=-\infty}^\infty W\skakko{M_n(\lambda+2\pi j)}$, and $$I_{ij}(\lambda):=\frac{1}{2\pi n}\skakko{\sum_{t=1}^{n}Y_i(t)e^{-\mathrm{i}\lambda t}}\skakko{\sum_{t=1}^{n}Y_j(t)e^{\mathrm{i}\lambda t}}.$$
Hence, it holds that
\begin{align}\label{eq_est_f}
{\hat{f}}_{ij}(\lambda)-{f}_{ij}(\lambda)=O_p\skakko{\sqrt\frac{M_n}{n}}
\quad \text{uniformly in $\lambda\in[-\pi,\pi]$}.
\end{align}
From \eqref{diff_f}, \eqref{eq_est_f}, and \citet[Lemma 3.4]{eichler08}, it can be seen that 
\begin{align}\nonumber
&\frac{n}{\sqrt {M_n}}\bigintssss_{-\pi}^\pi \left|\Phi_K\skakko{\hat {\bm{f}}(\lambda)}\right|^2{\rm d}\lambda-\frac{n}{\sqrt {M_n}}\bigintssss_{-\pi}^\pi \left|\Phi_K\skakko{\hat{\hat {\bm{f}}}(\lambda)}\right|^2{\rm d}\lambda\\\nonumber
=&\frac{n}{\sqrt {M_n}}
\int_{-\pi}^\pi {\rm vec}\skakko{\hat {\bm{f}}(\lambda)-{ {\bm{f}}}(\lambda)}^* 
\Gamma_{{\phi}}(\lambda)
{\rm vec}\skakko{\hat {\bm{f}}(\lambda)-{ {\bm{f}}}(\lambda)}\\\nonumber
&-
\skakko{\hat{\hat {\bm{f}}}(\lambda)-{{\bm{f}}}(\lambda)}^* 
\Gamma_{{\phi}}(\lambda)
{\rm vec}\skakko{\hat{\hat {\bm{f}}}(\lambda)-{ {\bm{f}}}(\lambda)}
{\rm d}\lambda+o_p\skakko{1}\\\nonumber
=&\frac{n}{\sqrt {M_n}}
\int_{-\pi}^\pi {\rm vec}\skakko{\hat {\bm{f}}(\lambda)-\hat{\hat {\bm{f}}}(\lambda)}^* 
\Gamma_{\Phi_K}(\lambda)
{\rm vec}\skakko{\hat {\bm{f}}(\lambda)-\hat{\hat {\bm{f}}}(\lambda)}\\\nonumber
&+
{\rm vec}\skakko{\hat {\bm{f}}(\lambda)- { {\bm{f}}}(\lambda)}^* 
\Gamma_{\Phi_K}(\lambda)
{\rm vec}\skakko{\hat {\bm{f}}(\lambda)-\hat{\hat {\bm{f}}}(\lambda)}\\\nonumber
&+
{\rm vec}\skakko{\hat {\bm{f}}(\lambda)-\hat{\hat {\bm{f}}}(\lambda)}^* 
\Gamma_{\Phi_K}(\lambda)
{\rm vec}\skakko{\hat {\bm{f}}(\lambda)- { {\bm{f}}}(\lambda)}{\rm d}\lambda
+o_p\skakko{1}
\\\label{diff_phi}
=&
O_p\skakko{\frac{M_n^{3/2}}{n}+ \frac{M_n}{n^{1/2}}}+o_p\skakko{1},
\end{align}
where
\begin{align*}
\Gamma_{ \Phi_K}(\lambda)
:= {\rm vec}
\left.
\skakko{\frac{\overline{\partial \Phi_K\skakko{{\bm Z}}}}{\partial {\bm Z}}}
{\rm vec}\skakko{\frac{{\partial \Phi_K\skakko{{\bm Z}}}}{\partial {\bm Z}}}^\top
\right|_{{\bm Z}={\bm{f}}(\lambda)}.
\end{align*}
Thus \eqref{diff_Tf} follows from \eqref{diif_mu} and \eqref{diff_phi}.

Second, we show the asymptotic normality of $T\skakko{\hat{\hat {\bm{f}}}(\lambda)}$. To this end, we check the conditions (i)--(iv) of \citet[Assumption 3.2]{eichler08}, that is, for the function $\Phi_K\skakko{\cdot}$ on the open set $\bm D:=\{{\bm Z}=(z_{ij})_{i,j=0,\ldots,K}\in \mathbb C^{(k+1)\times(k+1)};\|{\bm Z}\|<2L_D\}$ to $\mathbb C$ such that
\begin{align*}
\Phi_K\skakko{\bm Z}:=&-\sum_{i=1}^{K-1}{\rm det}\skakko{\overline{ {\bm{Z}}_{i,K}^\ddag}}z_{i0}+{\rm det}\skakko{{\bm{Z}}_{K-1}}z_{K0},
\end{align*}
where ${\bm{Z}}_{K-1}:=(z_{ij})_{i,j=1,\ldots,(K-1)}$ and 
\begin{align*}
{\bm{Z}}_{i,K}^\ddag:=
\begin{pmatrix}
z_{11}&\cdots&z_{1(i-1)}&z_{1K}&z_{1(i+1)}&\cdots&z_{1(K-1)}\\
\vdots&\vdots&\vdots&\vdots&\vdots&\vdots&\vdots\\
z_{(K-1)1}&\cdots&z_{(K-1)(i-1)}&z_{(K-1)K}&z_{(K-1)(i+1)}&\cdots&z_{(K-1)(K-1)}\\
\end{pmatrix},
\end{align*}
\begin{enumerate}
\item[(i)] The function $\Phi_K$ is holomorophic function.
\item[(ii)] The functions $\Phi_K(\bm{f}(\lambda))$ and 
$\left.{\rm vec}\skakko{\frac{{\partial \Phi_K\skakko{{\bm Z}}}}{\partial {\bm Z}}}\right|_{{\bm Z}={\bm{f}}(\lambda)}$ are Lipschitz continuous with respect to $\lambda\in[-\pi,\pi]$.
\item[(iii)] It holds that, for 
$B_{L_D,\lambda}:=\{{\bm Z}\in \mathbb C^{(k+1)\times(k+1)};\|{\bm Z}-{\bm{f}}(\lambda)\|<L_D\}$,
$$\sup_{\lambda\in[-\pi\pi]}\sup_{\bm Z\in B_{L_D,\lambda}} \left\|\Phi_K\skakko{{\bm Z}}\right\|<\infty.$$
\item[(iv)] The following holds true:
$$
\bigintssss_{-\pi}^\pi
\left\|\left.{\rm vec}\skakko{\frac{{\partial \Phi_K\skakko{{\bm Z}}}}{\partial {\bm Z}}}\right|_{{\bm Z}={\bm{f}}(\lambda)}\right\|>0.
$$
\end{enumerate}
From the definition of $\Phi_K$, (i) and (ii) are satisfied. 
The condition (iii) can be shown by
\begin{align*}
&\sup_{\lambda\in[-\pi,\pi]}\sup_{\bm Z\in B_{L_D,\lambda}} \left\|\Phi_K\skakko{{\bm Z}}\right\|\\
\leq&\sup_{\lambda\in[-\pi,\pi]}\sup_{\bm Z\in B_{L_D,\lambda}} \left\|\Phi_K\skakko{{\bm Z}}-\Phi_K\skakko{{\bm{f}}(\lambda)}\right\|
+\sup_{\lambda\in[-\pi,\pi]} \left\|\Phi_K\skakko{{\bm{f}}(\lambda)}\right\|\\
\leq&\sup_{\lambda\in[-\pi,\pi]}\sup_{\bm Z\in B_{L_D,\lambda}}L_{\Phi_K} \left\|{{\bm Z}}-{{\bm{f}}(\lambda)}\right\|
+\sup_{\lambda\in[-\pi,\pi]} \left\|\Phi_K\skakko{{\bm{f}}(\lambda)}\right\|\\
\leq&L_{\Phi_K}L_{D}+\sup_{\lambda\in[-\pi,\pi]} \left\|\Phi_K\skakko{{\bm{f}}(\lambda)}\right\|\\
<&\infty.
\end{align*}
Since $$\left.\skakko{\frac{{\partial \Phi_K\skakko{{\bm Z}}}}{\partial {z_{K0}}}}\right|_{{\bm Z}={\bm{f}}(\lambda)}=
{\rm det}\skakko{{\bm{f}}_{K-1}(\lambda)}\neq0,$$ (iv) holds.
\citet[Corollary 3.6.]{eichler08} yields the asymptotic normality of $T\skakko{\hat{\hat {\bm{f}}}(\lambda)}$. A lengthy algebra gives
the bias $\mu_K$ and asymptotic variance $\sigma_K^2$ as, for $K=1$, 
\begin{align*}
\mu_1:=&{\sqrt M_n}
\eta_{\omega,2}\int_{-\pi}^{\pi}f_{11}(\lambda)f_{00}(\lambda){\rm d}\lambda\quad\text{and }\sigma_1^2:=4\pi\eta_{\omega,4}\int_{-\pi}^{\pi}\left\| f_{11}(\lambda)f_{00}(\lambda)\right\|^{2}{\rm d}\lambda
\end{align*}
and, for $K\geq2$,
\begin{align*}
\mu_K:=&{\sqrt M_n}
\eta_{\omega,2}\int_{-\pi}^{\pi}(\det{\bm{f}_{K-1}(\lambda)})^{2}
\skakko{f_{00}(\lambda)-\overline{\hat{\bm{f}}_{K-1,0}^{\flat\top}}(\lambda)\bm{f}_{K-1}^{-1}(\lambda){\hat{\bm{f}}_{K-1,0}^\flat}(\lambda)}\\
&\quad\quad\quad\quad\quad\times
\skakko{f_{KK}(\lambda)-\overline{\hat{\bm{f}}_{K-1,K}^{\flat\top}}(\lambda)
\bm{f}_{K-1}^{-1}(\lambda)\hat{\bm{f}}_{K-1,K}^\flat(\lambda)}{\rm d}\lambda\\
\text{and }
\sigma_K^2:=&4\pi\eta_{\omega,4}\int_{-\pi}^{\pi}\left\| (\det{\bm{f}_{K-1}(\lambda)})^{2}
\skakko{f_{00}(\lambda)-\overline{\hat{\bm{f}}_{K-1,0}^{\flat\top}}(\lambda)\bm{f}_{K-1}^{-1}(\lambda){\hat{\bm{f}}_{K-1,0}^\flat}(\lambda)}\right.\\
&\quad\quad\quad\quad\quad\times\left.
\skakko{f_{KK}(\lambda)-\overline{\hat{\bm{f}}_{K-1,K}^{\flat\top}}(\lambda)
\bm{f}_{K-1}^{-1}(\lambda)\hat{\bm{f}}_{K-1,K}^\flat(\lambda)
}\right\|^{2}{\rm d}\lambda.
\end{align*}
\qed

\subsection{Proof of Theorem 3.2}
\citet[Theorem 5.1]{eichler08} yields, under the alternative $K_0$,
\begin{align*}
\frac{\sqrt {M_n}}{n}{T_n}
- 
\int_{-\pi}^\pi \left|\Phi_K\skakko{ {\bm{f}}(\lambda)}\right|^2{\rm d}\lambda
=o_{p}(1)\quad \text{as $n\to\infty$.}
\end{align*}
Then, we observe, under the alternative $K_0$,
\begin{align*}
{\rm P}\skakko{\frac{T_n}{\hat\sigma_K}\geq z_{\alpha}}
=&
{\rm P}\skakko{\frac{\sqrt {M_n}}{n}\frac{T_n}{\hat\sigma_K}\geq \frac{\sqrt {M_n}}{n}z_{\alpha}}\\
=&{\rm P}\skakko{\frac{1}{\sigma_K}\int_{-\pi}^\pi \left|\Phi_K\skakko{ {\bm{f}}(\lambda)}\right|^2{\rm d}\lambda
\geq 0}+o_p(1)\quad \text{as $n\to\infty$,}
\end{align*}
which shows the consistency of our test.\qed

\subsection{Proof of Theorem 3.3}
Since $\hat A_{ij}$ and $A_{ij}$ for $i,j(\leq i)$ are  functions of $\hat{\bm{f}}(\lambda)$ and ${\bm{f}}\skakko{\lambda}$, respectively, write $\hat A_{ij}$ and $A_{ij}$ as $\mathcal{A}_{ij}\skakko{\hat {\bm{f}}(\lambda)}$ and $\mathcal {A}_{ij}\skakko{{\bm{f}}(\lambda)}$, respectively.  
Note that $\mathcal {A}_{ij}$ is holomorphic, and thus, $\mathcal{A}_{ij}$ is Lipschitz continuous on the closed subset $\mkakko{{\bm Z}\in \mathbb C^{(K+1)\times(K+1)}:\left\|{\bm Z}-{\bm{f}}(\lambda)\right\|_2\leq \delta}$ for some $\delta>0$.
Therefore, we can show, for any $\epsilon>0$ and a Lipschitz constant $L$, 
\begin{align*}
&{\rm P}\skakko{\left|\hat a_{ij}(k)- a_{ij}(k)\right|>\epsilon}\\
=&
{\rm P}\skakko{\left|\frac{1}{2\pi}\int_{-\pi}^\pi \skakko{\mathcal{A}_{ij}\skakko{\hat {\bm{f}}(\lambda)}- \mathcal{A}_{ij}\skakko{ {\bm{f}}(\lambda)}}e^{-\mathrm{i}k\lambda}{\rm d}\lambda\right|>\epsilon}\\
\leq&
{\rm P}\skakko{\frac{L}{2\pi}\sup_{\lambda\in[-\pi,\pi]}
\left\|\hat {\bm{f}}(\lambda)-{\bm{f}}(\lambda)\right\|_2>\epsilon, 
\sup_{\lambda\in[-\pi,\pi]}\left\|\hat {\bm{f}}(\lambda)-{\bm{f}}(\lambda)\right\|\leq\delta
}\\
&+
{\rm P}\skakko{
\sup_{\lambda\in[-\pi,\pi]}\left\|\hat {\bm{f}}(\lambda)-{\bm{f}}(\lambda)\right\|_2>\delta
}
\end{align*}
which, in conjunction with $\sup_{\lambda\in[-\pi,\pi]}\left\|\hat {\bm{f}}(\lambda)-{\bm{f}}(\lambda)\right\|_2=o_p(1)$ (see, e.g., \citealp[Theorem 2.1]{robinson91}), tends to zero as $n\to\infty$. 

\qed

\end{appendix}



\begin{acks}[Acknowledgments]
We thank Brin Mathematics Research Center for supporting the first author's (Y.G.) research stay at the University of Maryland. We also thank Ms.\ Tong Lu for important assistance with data retrieval. Part of this research was conducted while the first author (Y.G.) was visiting the Brin Mathematics Research Center at the University of Maryland.
\end{acks}

\begin{funding}
This research was supported by JSPS Grant-in-Aid for Early-Career Scientists JP23K16851 (Y.G.) and NIH grant DP1 DA048968 (S.C.)
\end{funding}


\bibliographystyle{apalike}
\bibliography{ref}

\end{document}